\documentclass[12pt,a4paper]{article}
\usepackage{fullpage}
\usepackage {graphicx,subfigure, mathrsfs}
\usepackage{amsfonts,amsmath,amscd,amssymb,amsbsy,amsthm,amstext,amsopn, marvosym}
\usepackage{pstricks,pst-node,pst-text,pst-tree}
\usepackage[all]{xy}

\author{Tarig M. H. Abdelgadir\thanks{School of Mathematics, KIAS, 85 Hoegiro, Dongdaemun-gu, Seoul 130-722, Republic of Korea. $\text{\noindent Email: \texttt{tarig.m.h.abdelgadir@gmail.com}}$}}
\title {Quivers of sections on toric orbifolds}
\date {\today}

\newcommand{\ee}{\mathbf{e}}

\newcommand{\kk}{\ensuremath{\Bbbk}}
\newcommand{\pic}{\textup{pic}}
\newcommand{\Pic}{\textup{Pic}}
\renewcommand{\ker}{\textup{ker}}
\newcommand{\Spec}{\textup{Spec}}
\newcommand{\rank}{\textup{rank}}
\renewcommand{\div}{\textup{div}}
\newcommand{\Aut}{\textup{Aut}}
\newcommand{\Hom}{\textup{Hom}}
\newcommand{\Wt}{\textup{Wt}}
\newcommand{\GL}{\textup{GL}}
\newcommand{\SL}{\textup{SL}}
\newcommand{\PGL}{\textup{PGL}}
\newcommand{\inc}{\textup{inc}}
\newcommand{\ZZ}{\mathbb{Z}}
\newcommand{\RR}{\mathcal{R}}
\newcommand{\QQ}{\mathbb{Q}}
\newcommand{\BB}{\mathfrak{B}}

\newcommand{\VV}{\mathbb{V}}
\renewcommand{\AA}{\mathbb{A}}
\newcommand{\PP}{\mathbb{P}}
\newcommand{\WW}{\mathcal{W}}
\newcommand{\MM}{\mathcal{M}}

\newcommand{\LL}{\mathscr{L}}

\renewcommand{\b}{b}
\newcommand{\OO}{\mathcal{O}}
\newcommand{\XX}{\mathcal{X}}

\newcommand{\bsf}{base-point free }
\newcommand{\rep}{representation }
\newcommand{\reps}{representations }
\newcommand{\cms}{coarse moduli space }
\newcommand{\NN}{\mathbb{N}}
\renewcommand{\iff}{if and only if }

\theoremstyle {plain}

\newtheorem {thm}{Theorem}[section]
\newtheorem* {thm3.8}{Theorem 3.5}
\newtheorem*{thm4.12}{Theorem 4.13}
\newtheorem*{thm4.21}{Theorem 4.18}
\newtheorem*{thm5.4}{Theorem 5.4}
\newtheorem*{coro5.5}{Corollary 5.5}
\newtheorem {lemma}[thm]{Lemma}
\newtheorem {prop}[thm]{Proposition}
\newtheorem {coro}[thm]{Corollary}

\theoremstyle {definition}

\newtheorem {ex}[thm]{Example}
\newtheorem* {rmk4.7}{Remark 4.8}
\newtheorem {defn}[thm]{Definition}
\newtheorem {rmk}[thm]{Remark}

\newtheorem {conventions}[thm]{Conventions}

\begin{document}
\maketitle	

\abstract{Starting from a collection of line bundles on a projective toric orbifold $\XX$, we introduce a stacky analogue of the classical linear series. Our first main result extends work of King by building moduli stacks of refined representations of labelled quivers. We associate one such stack to any collection of line bundles on $\XX$ to obtain our notion of a stacky linear series; as in the classical case, $\XX$ maps to the ambient stack by evaluating sections of line bundles in the collection. As a further application, we describe a finite sequence of GIT wall crossings between $[\AA^n/G]$ and $G$-Hilb$(\AA^n)$ for a finite abelian subgroup $G \subset \SL(n,\kk)$ where $n \leq 3$.}

\section{Introduction} 

This article gives an analogue of the linear series construction for projective toric orbifolds. We start with a finite collection of line bundles $\LL$ and use quiver of sections, as defined by Craw-Smith \cite{Craw-Smith}, to package efficiently  the sections of line bundles in $\LL$. For our ambient space, we introduce moduli stacks $\MM_\theta(Q,\div)$ of quiver representation-like objects and produce rational maps $\psi_\theta: \XX \dashrightarrow \MM_\theta(Q, \div)$. As a further application, we apply our construction to orbifolds $[\AA^n/G]$, where $G$ is a finite abelian subgroup of $\GL(n,\kk)$, recovering the orbifold from the corresponding McKay quiver. We then adapt of our construction to describe a finite sequence of GIT wall crossings between $[\AA^n/G]$ and $G$-Hilb$(\AA^n)$ for $G \subset \SL(n,\kk)$ where $n \leq 3$.

Motivated by Olsson-Starr \cite{Olsson-Starr} among others, Kresch \cite{Kresch} introduced the notion of a projective Deligne-Mumford stack. A Deligne-Mumford stack $\XX$ is said to be projective if it has a projective coarse moduli space $X$ and a generating sheaf. One may think of a generating (or $\pi$-very ample) sheaf as a very ample sheaf relative to the morphism $\pi: \XX \rightarrow X$, or more loosely, as a sheaf that allows one to lift the projectivity of $X$ to $\XX$. Despite the importance of linear series in the theory of algebraic varieties, an appropriate analogue for algebraic stacks does not exist. The naive extension to algebraic stacks is not satisfactory. Indeed, requiring a line bundle to be both very ample on the \cms and $\pi$-very ample is too restrictive; the stack must be an algebraic space, forcing all stabilizers to be trivial. 

One could sidestep this issue by considering sections of more than one line bundle. In the case where $\XX$ has cyclic stabilizers, the approach adopted by Abramovich-Hassett \cite{AbramovichHassett} uses sections of tensor powers of a single line bundle $L$ to produce closed immersions $\XX \hookrightarrow \PP(\bigoplus_{j=n} ^m \Gamma(\XX, L^{\otimes j}))$ for some $n,m \in \NN$. For $\XX$ a toric orbifold, this article gives an alternative stacky analogue to the linear series construction that generalizes the Abramovich-Hassett construction. Our construction puts no constraints on the stabilizers. Moreover, the efficiency of the quiver of sections allows for ambient spaces of dimensions smaller than those appearing in the Abramovich-Hassett construction (cf. Example \ref{p112ex} and Example \ref{AHex1}).

Our construction is not limited to projective stacks. In fact for $G \subset \GL(n, \kk)$ a finite abelian group, it recovers the stack quotient $[\AA^n / G]$ from the McKay quiver. Under some constraints on $G$, Craw-Ishii \cite{CrawIshii} show that the McKay quiver allows us to move between projective crepant resolutions of $\AA^n/G$ by a finite sequence of wall crossings. When $G \subset \SL(n,\kk)$, one may think of the stack $[\AA^n/ G]$ as a noncommutative crepant resolution of $\AA^n / G$. This is because its coordinate ring, as defined by Chan-Ingalls \cite{ChanIngalls}, is a noncommutative crepant resolution of $\AA^n/G$ in the sense of Van den Bergh \cite{NoncommCrepant}. Therefore it is natural to ask whether one can introduce a quiver theoretic construction that allows us to move between a crepant resolution of  $\AA^n/G$, say Nakamura's $G$-Hilb$(\AA^n)$, and the stack $[\AA^n/ G]$ by crossing finitely many walls. A slight adaptation of our construction gives an affirmative answer to this question, putting $G$-Hilb$(\AA^n)$ and $[\AA^n/G]$ on the same footing.

We now summarize the contents of the paper in more detail. Motivated by the natural labelling of the quiver of sections on a toric variety by torus-invariant divisors, we define the notion of a labelled quiver. A {\em labelled quiver} is a quiver $Q$ along with a map of sets $l:Q_1 \rightarrow \ZZ^d$. Naively, one wishes to define a `representation' of a labelled quiver as a representation of the underlying quiver for which any two paths with the same label are represented by the `same' linear map. Forcing this on linear maps representing two paths that have the same labels but don't share the same head and tail introduces equations that are not homogenous, with respect to the change of basis action, on the representation space. We bypass this issue by introducing new `homogenizing' parameters to the representation space, that homogenize every equation of paths induced by the labels. A {\em refined representation} of a labelled quiver $(Q,l)$ is a representation of the underlying quiver, together with a choice of nonzero homogenizing parameters. Given a refined representation $W$ and weight $\theta \in K_0(\kk Q \text{-mod})^\vee$, we define a notion of $\theta$-stability on $W$ following King \cite{King}. For a weight $\theta \in K_0(\kk Q \text{-mod})^\vee$ defined by a character $\chi_\theta$ of $\PGL(\alpha)$ (the group acting faithfully on the refined representation space), the main result of Section \ref{Refine} relates GIT $\chi_\theta$-stability to $\theta$-stability.
\begin{thm3.8}
Let $\chi_\theta$ be a character of $\GL(\alpha)$ and $\theta$ the corresponding element of $K_0(\kk Q \text{-mod})^\vee$. A refined quiver representation $W$  is $\theta$-semistable (resp. $\theta$-stable) if and only if the corresponding point in $\RR(Q, l, \alpha)$ is $\chi_\theta$-semistable (resp. $\chi_\theta$-stable) with respect to action of $\GL(\alpha)$.
\end{thm3.8}
\noindent This allows us to introduce families of $\theta$-semistable refined representations, which in turn enables us to define moduli stacks $\MM_\theta(Q, l, \alpha)$ of refined representations, given some dimension vector $\alpha$. The stacks $\MM_\theta(Q, l, \alpha)$ form the ambient stacks in our construction.

Now take $\XX$ to be a projective toric orbifold. For a given collection of line bundles $\LL$ on $\XX$, we use techniques very similar to those in \cite{Craw-Smith} to define a labelled quiver of sections $(Q,\div)$ and give a rational map $\psi_\theta: \XX \dashrightarrow \MM_\theta(Q,\div, \alpha)$ where $\alpha = (1,\ldots,1)$. As in the classical linear series case, when $\psi_\theta$ is a morphism the tautological line bundles on $\MM_\theta(Q,\div, \alpha)$ pull-back to recover the collection $\LL$.  Checking whether or not there exists a stability condition $\theta$ for which $\psi_\theta$ is a morphism can be tedious, hence we introduce a sufficient condition that is straightforward to check. We also explicitly describe the image of $\psi_\theta$ and address the question of representability of the morphism $\psi_\theta$. Let $\LL_{\textup{bpf}}$ denote the collection of line bundles \[\{L_i^\vee \otimes L_j\, |\, L_i, L_j \in \LL \text{ and } L_i^\vee \otimes L_j \text{ is base-point free}\}\] we show the following,

\begin{thm4.12}
If $\rank(\ZZ \LL) = \rank(\ZZ \LL_\textup{bpf})$ then $\LL$ is base-point free, i.e.\ there exists a generic stability condition $\theta$ such that $\psi_\theta: \XX \rightarrow \MM_\theta(Q,\div)$ is a morphism. 
\end{thm4.12}

\begin{thm4.21}
A morphism $\psi_{\theta}$ is representable \iff $\bigoplus_{j=1}^rL_j$ is $\pi$-ample.
\end{thm4.21}

For $\XX$ a toric orbifold, the Abramovich-Hassett construction may be recovered.
\begin{rmk4.7}
Given a polarizing line bundle $L$ on $\XX$,  the Abramovich-Hassett construction for $n=0$, is recovered by applying our machinery to the collection  \[\LL = (\OO_\XX, L, L \otimes L^{\otimes 2}, \ldots, L^{\otimes m(m+1)/2})\]  and if necessary, working with an `incomplete' quiver of sections. An incomplete quiver of sections is a quiver of sections where not all torus-invariant sections contribute to paths in the quiver, analogous to an incomplete linear series.
\end{rmk4.7}

We then apply this technology to the McKay quiver associated to a finite abelian subgroup of $\GL(n, \kk)$. After showing that every refined representation of the labelled McKay quiver $(Q, \div)$ is $\theta$-stable and deducing $\psi_\theta$ is a morphism for any given stability condition $\theta$, we show that $\psi_\theta: \AA^n/G \rightarrow \MM_\theta(Q, \div)$ is a closed immersion. We tweak the GIT construction of $\MM_\theta(Q, \div)$ by allowing the homogenizing parameters to be zero and examine a substack cut out by an ideal defined naturally from the labels of $Q$. By studying the GIT chamber decomposition, we observe that certain chambers define semistable loci in which every homogenizing parameter is nonzero, enabling us to recover the stack $[\AA^n/G]$. We also show that in the semistable locus of a second chamber the homogenizing variables are completely determined by the variables corresponding to the arrows and are therefore redundant. This recovers the Craw-Maclagan-Thomas \cite{CMT1} construction of the coherent component Hilb$^G(\AA^n)$  of Nakamura's $G$-Hilbert scheme. Using the results of Ito-Nakamura \cite{ItoNakamura} and Nakamura \cite{Nakamura} we have the following results.
\begin{thm5.4}
For finite abelian $G \subset \GL(n,\kk)$, there exists generic stability conditions $\chi_{\theta_1}, \chi_{\theta_2} \in \PGL(\alpha)^\vee$, such that \[[\AA^n/G] \cong [\VV(I_{\LL, \BB})^{ss}_{\theta_1} / \PGL(\alpha)] \,\text{ and }\, \textup{Hilb}^G(\AA^n) \cong [\VV(I_{\LL, \BB})^{ss}_{\theta_2} / \PGL(\alpha)].\]
\end{thm5.4}

\begin{coro5.5}
For $n \leq 3$ and finite abelian $G\subset \SL(n, \kk)$, there exists generic stability conditions $\chi_{\theta_1}, \chi_{\theta_2} \in \PGL(\alpha)^\vee$, such that \[ [\AA^n/G] \cong [\VV(I_{\LL, \BB})^{ss}_{\theta_1} / \PGL(\alpha)] \,\text{ and }\, G\textup{-Hilb}(\AA^n) \cong [\VV(I_{\LL, \BB})^{ss}_{\theta_2} / \PGL(\alpha)].\]
\end{coro5.5}

The paper is organized as follows. We assemble some background material on quivers and toric orbifolds in Section \ref{Back}. In Section \ref{Refine}, we define the ambient stacks in our construction. We define our stacky analogue of the classical linear series construction in Section \ref{Quiver}. Finally in Section \ref{McKay}, we apply the machinery from Section \ref{Quiver} to the McKay quiver.

\subsubsection*{Conventions and Notation} The symbol $\kk$ will be reserved for an algebraically closed field of characteristic 0. All objects and maps are defined over $\kk$ unless stated. The symbol $\NN$ will be reserved for the nonnegative integers. Our semigroups contain units. For a finite set $C$ we use $\ZZ C$ to denote the free abelian group generated by $C$ and $\NN C$ to be free the abelian semigroup generated by $C$. For an abelian group $G$ we write $G_\QQ:= G \otimes_\ZZ \QQ$. For us, a geometric point of a stack $\XX$ is a morphism $\Spec(\kk) \rightarrow \XX$. We use $\PP(w_1, \ldots, w_n)$ to denote the weighted projective stack with weights $w_1,\ldots,w_n$. 

\subsubsection*{Acknowledgments}
Firstly, I would like to express my gratitude to my Ph.D. supervisor Alastair Craw for introducing me to this problem and for the countless discussions that shaped this project. I would like to thank Alastair King for many stimulating discussions and for his ideas on generalizing refined representations beyond the case $\alpha = (1,\ldots,1)$. I am indebted to Barbara Fantechi and {\'E}tienne Mann for teaching me about stacks and answering numerous stacky questions. I would also like to thank Ewan Morrison, Alexander Quintero V\'{e}lez,  Dorothy Winn and the rest of office 522 for all the helpful discussions we had. Many thanks are due to the referee and my thesis examiners for useful comments and suggestions. I also acknowledge the financial support of the University of Glasgow.

\section{Background}\label{Back}

\subsection{Quivers}\label{Quivers}
Most of this subsection is borrowed from Section 2.2 of \cite{Craw-Smith} and is included here for completeness.

A quiver $Q$ is specified by two finite sets $Q_0$ and $Q_1$, whose elements are called vertices and arrows, together with two maps $h, t \colon Q_1 \rightarrow Q_0$ indicating the vertices at the head and tail of each arrow.  We assume through out that $Q$ is connected. A nontrivial path in $Q$ is a sequence of arrows $p = a_1 \dotsb a_m$ with $h(a_{k}) = t(a_{k+1})$ for $1 \leq k < m$.  We set $t(p) = t(a_{1})$ and $h(p)= h(a_m)$.  Each $i \in Q_0$ gives a trivial path $e_i$ where $t(e_i) = h(e_i) = i$.  The path algebra $\kk Q$ is the $\kk$-algebra whose underlying $\kk$-vector space has a basis consisting of paths in $Q$; the product of two basis elements equals the basis element defined by concatenation of the paths if possible or zero otherwise.  A cycle is a path $p$ in which $t(p) = h(p)$. A quiver is acyclic if it contains no cycles.  A vertex is a source if it is not the head of any arrow and a quiver is rooted if it has a unique source. 

The vertex space $\ZZ^{Q_0}$ is the free abelian group generated by the vertices and the arrow space $\ZZ^{Q_1}$ is the free abelian group generated by the arrows. We write $\NN^{Q_0}$ and $\NN^{Q_1}$ for the subsemigroups generated by the basis elements of $\ZZ^{Q_0}$ and $\ZZ^{Q_1}$. The incidence map $\inc: \ZZ^{Q_1} \rightarrow \ZZ^{Q_0}$ is defined by $\inc(\ee_a) = \ee_{h(a)} - \ee_{t(a)}$. The weight lattice $\Wt(Q)$ is the image of $\inc$, that is the sublattice give by elements $\theta = \sum_{i\in Q_0} \theta_i \ee_i \in \ZZ^{Q_0}$ for which $\sum_{i \in Q_0} \theta_i = 0$.

A representation $\overline{W} = (W_i, w_a)$ of $Q$ consists of a vector space $W_i$ for each $i \in Q_0$ and a linear map $w_a \colon W_{t(a)} \rightarrow W_{h(a)}$ for each $a \in Q_1$. The dimension vector of $\overline{W}$ is the integer vector $(\dim W_{i}) \in \NN^{Q_0}$.  A map between representations $\overline{W} = (W_i, w_a)$ and $\overline{W}' = (W_i', w_a')$ is a family $\xi_{i} \colon W_i^{\,} \rightarrow W_i'$ for $i \in Q_0$ of linear maps that are compatible with the structure maps, that is $w_a' \circ\xi_{t(a)} = \xi_{h(a)} \circ w_a$ for all $a \in Q_1$.  With composition defined componentwise, we obtain the abelian category of representations of $Q$ denoted rep$_\kk(Q)$. This category is equivalent to the category $\kk Q$-mod of finitely generated left modules over the path algebra.

Each $\theta \in \Wt(Q)$ defines a stability notion for representations of $Q$.  A representation $\overline{W}$ is $\theta$-semistable if, for every proper, nonzero subrepresentation $\overline{W}' \subset \overline{W}$, we have $\sum_{i \in Q_0} \theta_i \cdot \dim(W_i') \geq 0$.  The notion of $\theta$-stability is obtained by replacing $\geq$ with $>$. For a given dimension vector $\alpha \in \NN^{Q_0}$, a family of $\theta$-semistable quiver representations over a connected scheme $S$ is a collection of rank $\alpha_i$ locally free sheaves $\WW_i$ together with morphisms $\WW_{t(a)} \rightarrow \WW_{h(a)}$ for every $a\in Q_1$. When every $\theta$-semistable representation is $\theta$-stable and the dimension vector is primitive this moduli problem is representable by a scheme $\MM_\theta(Q, \alpha)$, see Proposition 5.3 in \cite{King}. 

\subsection{Toric orbifolds}

Toric Deligne-Mumford stacks were first introduced by Borisov-Chen-Smith in \cite{BCS}, later Fantechi-Mann-Nironi \cite{FMN} gave an equivalent definition analogous to the classical definition of a toric variety. In this article we will only be concerned with toric orbifolds, that is smooth toric Deligne-Mumford stacks whose generic stabilizer is trivial or equivalently stacks whose dense Deligne-Mumford torus is just an algebraic torus $T$. We begin by introducing the Fantechi-Mann-Nironi approach then discuss the Borisov-Chen-Smith approach.

A toric orbifold is a smooth separated Deligne-Mumford stack $\XX$ together with an open immersion $\iota: T \hookrightarrow \XX$ with dense image such that the action of $T$ on itself extends to an action of $T$ on $\XX$. Toric orbifolds can also be defined using stacky fans. A stacky fan is a triple $\mathbf{\Sigma} := (N, \Sigma, \beta)$, where $N$ is a finitely  generated free abelian group, $\Sigma$ is a rational simplicial fan in $N_\QQ$ with $d$ rays that span $N_\QQ$, denoted $\rho_1, \ldots, \rho_d \in \Sigma(1)$, and $\beta: \ZZ^d \rightarrow N$ is a morphism of groups for which $\beta(\ee_i) \otimes 1$ is on the ray $\rho_i \in N_\QQ$. The toric orbifold associated to a stacky fan is constructed as follows. Let $\ZZ^{\Sigma(1)}:= (\ZZ^d)^\vee$ and consider the exact sequence
\begin{equation}\label{Coxsq}
\xymatrix@C=1.3cm{N^\vee \ar[r]^-{\beta^\vee} & \ZZ^{\Sigma(1)} \ar[r]^-{\deg} & \text{Coker}(\beta^\vee) \ar[r]& 0.}
\end{equation}
The the group $T = \Hom(\text{Coker}(\beta^\vee), \kk^\times)$ has a natural action on $\AA^{\Sigma(1)}$ induced by the inclusion $T \subset (\kk^\times)^{\Sigma(1)}$. For a cone $\sigma \in \Sigma$, $\widehat{\sigma}$ is the set of one-dimentional cones in $\Sigma$ not contained in $\sigma$ and $x^{\widehat{\sigma}} = \prod_{\rho \in \widehat{\sigma}} x_\rho$. The Cox unstable locus is 
the defined 
\begin{equation}\label{Coxunstable}
B_\XX:= \Big\langle x^{\widehat{\sigma}} \in \kk[ x_{\rho} \, | \, \rho \in \Sigma(1)] \,\,\Big | \,\,\sigma \in \Sigma\Big\rangle.
\end{equation}
The stack $\XX_\mathbf{\Sigma}$ associated to the stacky fan is \[\Bigg[ \frac{\AA^{\Sigma(1)} \setminus \VV(B_\XX)}{T}\Bigg].\] The group of line bundles $\Pic(\XX_{\mathbf{\Sigma}})$ is given by $\Hom(T, \kk^*) \cong \text{Coker}(\beta^\vee)$ and the group of torus-invariant divisors of $\XX_\mathbf{\Sigma}$ is given by $\ZZ^{\Sigma(1)}$. Given any toric orbifold $\XX$ there exists a stacky fan $\mathbf{\Sigma}$ such that $\XX_\mathbf{\Sigma} \cong \XX$, see Theorem 7.23 in \cite{FMN}. A toric orbifold is projective, in the sense of Kresch \cite{Kresch}, if its coarse moduli space is projective (cf. Corollary 5.4, \cite{Kresch}).

\section{Moduli of refined quiver representations}\label{Refine}
The goal of this section is to define an ambient space for our stacky analogue of the linear series construction. A quiver of sections of a toric orbifold is naturally labelled by torus-invariant divisors (see Section \ref{Quiver}). Motivated by this, we define the notion of a labelled quiver of which our quivers of sections are examples. Then we define the notion of a refined quiver representation of a labelled quiver and study moduli stacks of the aforementioned objects.

\begin{defn}
A {\em labelled quiver} $(Q, l)$ is a connected finite quiver $Q$ along with a free abelian group $\ZZ^d$ for some $d\in \NN$ and a map of sets $l: Q_1 \rightarrow \ZZ^d$.
\end{defn}

Abusing notation we use $l$ to denote the {\em labelling map} $\ZZ^{Q_1} \rightarrow \ZZ^d$ generated by $l$. Let $R$ denote the image of $\ker(l)$ under the incidence map and consider the following commutative diagram,
\begin{equation}
\begin{split}
\xymatrix{
\text{ker}(l) \ar[r] \ar[d] & R:= \text{inc(ker}(l)) \ar[d]\\
\ZZ^{Q_1} \ar[r]^-{\text{inc}} \ar[d]_{l} & \Wt(Q)\\
\ZZ^{d}.}
\end{split}
\end{equation}
Given a quiver $Q$, a representation $(W_i, w_a)$ of $Q$ and an element $\b = \sum_{i \in Q_0} \b_i \ee_i \in \ZZ^{Q_0}$ define
\begin{equation*}
\textup{det}_{\b}\, W:= \bigotimes_{i\in Q_0}(\det W_i)^{\otimes \b_i}.
\end{equation*}
Here we use the convention $W_i ^{\otimes -1} := W_i^{\vee}$.

Pick a basis $\BB$ of $R$.

\begin{defn} \label{rerep}
A {\em refined representation $W$ of a labelled quiver} $(Q,l)$ consists of a finite dimensional representation $\overline{W}:=(W_i, w_a)$ of $Q$ together with an isomorphism $f_\b: \kk \rightarrow \text{det}_b\,W$ for every $\b \in \BB$. The {\em dimension vector} of an $R$-refined representation $W= (W_i, w_a, f_\b)$ is the integer vector $(\text{dim}(W_i))_{i \in Q_0}$.

We say that two refined \reps $W = (W_i, w_a, f_\b),\, W' = (W_i', w_a', f'_\b)$ are {\em isomorphic} if there exist isomorphisms of vector spaces $\gamma_i: W_i \rightarrow W_i'$ for every vertex $i \in Q_0$ such that $\gamma_{t(a)}^{-1} \circ w_a \circ \gamma_{h(a)} = w_a'$ for all $a\in Q_1$and $f_\b \circ \gamma_\b =f'_\b$ for all $\b \in \BB$ where $\gamma_b: \text{det}_b\,W \rightarrow \text{det}_b\,W'$ is the isomorphism induced by the isomorphisms $\gamma_i$. 
\end{defn}

\begin{rmk}
\begin{enumerate} 
\item [i)] The independence of the choice of basis $\BB$ will be addressed in Remark \ref{indep}.
\item [ii)] Refined representations and their moduli maybe defined without appealing to a labelling map $l$; the crucial ingredient is the subgroup $R \subset \Wt(Q)$. Given a quiver $Q$ and an arbitrary subgroup $K \subset \Wt(Q)$ with a choice of basis $\BB_K$, one may define {\em $K$-refined representations of $Q$} to be finite dimensional representation $\overline{W}$ of $Q$ together isomorphisms $f_\b: \kk \rightarrow \text{det}_{\b}\,W$ for every $\b \in \BB_K$. All the definitions and results in this section may be lifted to this setting. With the immediate applications in mind, we restrict ourselves to subgroups $R$ arising from a labelling map $l$.
\end{enumerate}
\end{rmk}

For $i\in {Q_0}$, let $W_i$ be a vector space of dimension $\alpha_i$ and $\alpha:= (\alpha_i) \in \NN^{Q_0}$. Let $\kk[z_\b\,|\, \b \in \BB]$ denote the coordinate ring of the vector space $\bigoplus_{\b \in \BB}\text{det}_{\b}\, W$. The isomorphism classes of refined representations of $(Q, l)$ are in one-to-one correspondence with the orbits in the refined representation space 
\begin{equation*}
\RR(Q, l, \alpha) := \bigg( \bigoplus_{a \in Q_1} \Hom(W_{t(a)}, W_{h(a)}) \oplus \bigoplus_{\b \in \BB}\text{det}_{\b}\, W \bigg) \setminus \VV\Big(\prod_{\b \in \BB} z_\b\Big)
\end{equation*}
of the symmetry group
\begin{equation*}
\GL(\alpha) := \prod_{i \in Q_0} \GL(W_i)
\end{equation*}
under the change of basis action. Note that $\GL(\alpha)$ contains the diagonal one-parameter subgroup $\Delta = \{(\lambda \cdot 1,\ldots,\lambda\cdot 1): \lambda \in \kk^\times\}$ acting trivially and define $\PGL(\alpha) := \GL(\alpha) / \Delta$. 

We note that the characters of $\GL(\alpha)$ are given by
\begin{equation*}
\chi_\theta(g) = \prod_{i \in Q_0} \det(g_i)^{\theta_i}
\end{equation*}
for $\theta= \sum_i \theta_i \ee_i \in \ZZ^{Q_0}$ and that every character of $\GL(\alpha)$ is of the form $\chi_\theta$ for some $\theta \in \ZZ^{Q_0}$. As the diagonal $\Delta \subset \GL(\alpha)$ acts trivially on $\RR(Q, \div, \alpha)$ we are interested in characters $\chi_\theta$ that satisfy $\sum_i \theta_i \alpha_i =0$.

It is convenient to identify $\ZZ^{Q_0}$, and hence $\GL(\alpha)^\vee$, with a subgroup of the Grothendieck group $K_0(\kk Q \text{-mod})$ as follows. Let $\overline{W}=(W_i, w_a)$ be a representation of $Q$, implicitly using $\kk Q \text{-mod} \cong \text{rep}_\kk(Q)$, set $\theta(\overline{W}) = \sum_i \theta_i \text{dim }W_i$, and observe that this is additive on short exact sequences.

We introduce some notation before the next definition. Let $M$ a module of some ring $R$ and $M_\bullet$ be a proper filtration of $M$ (that is, a filtration where at least one term is a nonzero proper submodule of $M$) given by
\begin{equation*}
0 \subsetneq M_1 \subset \ldots \subset M_{n-1} \subsetneq M_n=M.
\end{equation*}
For $\theta \in K_0(\text{mod-}R)^\vee$ define $\theta(M_\bullet)= \sum_{j=1}^{n-1} \theta(M_j)$.

\begin{defn} \label{stab}
Let $\theta \in K_0(\kk Q \text{-mod})^\vee$. A refined representation $W$ is {\em $\theta$-semistable} if $\theta(W) = 0$ and $\theta(W_\bullet)\geq 0$ for every proper filtration $W_\bullet$ of the $\kk Q$-module $\overline{W}$ that satisfies $\b(W_\bullet) = 0$ for every $\b \in \BB$.
The notion of {\em $\theta$-stability} is obtained by replacing $\geq$ with $>$.
\end{defn}

We introduced the notion of $\theta$-semistability to be able to make sense of families of refined representations and use the term `moduli stack'. In practice, we are interested primarily in functions $\theta \in K_0(\kk Q \text{-mod})^\vee$ coming from characters of $\GL(\alpha)$. The next proposition allows us to check $\theta$-semistability by checking GIT semistability of $\chi_\theta$.

\begin{thm}\label{thetaGIT}
Let $\chi_\theta$ be a character of $\GL(\alpha)$ and $\theta$ the corresponding element of $K_0(\kk Q \text{-mod})^\vee$. A refined quiver representation $W$  is $\theta$-semistable (resp. $\theta$-stable) if and only if the corresponding point in $\RR(Q, l, \alpha)$ is $\chi_\theta$-semistable (resp. $\chi_\theta$-stable) with respect to action of $\GL(\alpha)$.
\end{thm}

\begin{proof}
We begin by pinning down the one-parameter subgroups $\lambda$ of $\GL(\alpha)$ for which $\lim_{t \rightarrow 0} (\lambda(t)\cdot W)$ exists. Write $\RR(Q, l, \alpha) \cong \RR(Q, \alpha) \times (\kk^\times)^\BB$ and $\pi_1, \pi_2$ for the first and second projection respectively. The limit $\lim_{t \rightarrow 0} (\lambda(t)\cdot W)$ exists \iff $\lim_{t \rightarrow 0} (\lambda(t)\cdot \pi_1(W))$ and $\lim_{t \rightarrow 0} (\lambda(t)\cdot \pi_2(W))$ exist. By the discussion preceding Proposition 3.1 of \cite{King}, $\lim_{t \rightarrow 0} (\lambda(t)\cdot \pi_1(W))$ exists if and only if $\lambda$ defines a $\ZZ$-filtration, $W_\bullet$, of the $\kk Q$-module $\pi_1(W) = \overline{W}$,
\begin{equation*}
\ldots \subset W_{n-1} \subset W_n \subset W_{n+1} \subset \ldots
\end{equation*}
for which $W_n = 0$ for $n \ll 0$ and $W_n = \overline{W}$ for $n \gg0$. Now consider $\lim_{t \rightarrow 0} (\lambda(t)\cdot \pi_2(W))$. The one-parameter subgroup  $\lambda$ defines a $\ZZ$-grading on the coordinate ring $\kk[z_{\b}, z_\b^{-1} \, | \, \b \in \BB]$ of $(\kk^\times)^\BB$. The limit  $\lim_{t \rightarrow 0} (\lambda(t)\cdot \pi_2(W))$ exists if and only if the variables $z_\b$ and $z_\b^{-1}$ are simultaneously non-negatively graded. Notice that this holds precisely when they are zero graded, that is when $\langle \chi_\b, \lambda\rangle =0$ for every $\b \in \BB$. Therefore, for $\lambda$ and $W$ as above, $\lim_{t \rightarrow 0} (\lambda(t)\cdot W)$ exists \iff $\lambda$ gives a $\ZZ$-filtration $(W_n)_{n \in \ZZ}$ of the quiver \rep$\pi_1(W) =\overline{W}$ and $\langle \chi_\b, \lambda\rangle =0$, for every $\b \in \BB$. 

Now assume $W$ is $\theta$-semistable. Take $\lambda$ to be a one-parameter subgroup for which the limit $\lim_{t \rightarrow 0} (\lambda(t)\cdot W)$ exists. By the discussion preceding Proposition 3.1 of \cite{King}, one may associate a filtration $W_\bullet$ to $\lambda$ such that $\langle \chi_\theta, \lambda\rangle = \theta(W_\bullet)$. By assumption, this implies $\langle \chi_\b, \lambda\rangle = \b(W_\bullet) = 0$ for all $\b \in \BB$. Since $W$ is $\theta$-semistable we have $\langle \chi_\theta, \lambda\rangle = \theta(W_\bullet) \geq 0$. GIT semistability of $W$ then follows from Mumford's numerical criterion, see Proposition 2.5 of \cite{King}.

Next assume $W \in \RR(Q, l, \alpha)$ is $\chi_\theta$-semistable. By the fact that $\Delta$ acts trivially we have $\langle \chi_\theta, \Delta\rangle = \theta(W) = 0$. Let $W_\bullet$ be a proper filtration satisfying the conditions of Definition \ref{stab}. By the discussion preceding Proposition 3.1 of \cite{King}, there exists a one-parameter subgroup $\lambda$ for which the associated filtration is $W_\bullet$. By assumption we have that $\b(W_\bullet) = \langle \chi_\b, \lambda\rangle = 0$ for every $\b \in \BB$, so $\lim_{t \rightarrow 0} (\lambda(t)\cdot W)$ exists. Mumford's numerical criterion gives $\theta(W_\bullet) =  \langle \chi_\theta, \lambda\rangle \geq 0$, as required.
\end{proof}

\begin{defn}\label{moduli}
For $\chi_\theta \in \PGL(\alpha)^\vee \subset \GL(\alpha)^\vee$, let $\RR(Q, l, \alpha)_\theta^{ss}$ denote the open subscheme of $\RR(Q, l, \alpha)$ parametrizing the $\theta$-semistable refined representation. The {\em moduli stack of $\theta$-semistable refined representations } is the stack quotient 
\begin{equation*}
\MM_\theta(Q, l, \alpha):= [\RR(Q, l, \alpha)_\theta^{ss} / \PGL(\alpha)].
\end{equation*}
\end{defn}

\begin{rmk} \label{indep}
The definition of $\MM_\theta(Q, l, \alpha)$ depends a priori on a choice of basis $\BB$ of $R$. However, any alternative basis $\BB'$ gives an isomorphic stack. Indeed, given $W=(W_i, w_a, f_\b) \in \RR(Q, l, \alpha)$ write $\b' = n_1 \b_1 + \cdots + n_m \b_m$ and \[f_{\b'} :\kk \cong \kk^{\otimes n_1}\otimes \cdots \otimes \kk^{\otimes n_m}  \rightarrow (\text{det}_{b_1} \,W)^{\otimes n_1}\otimes \cdots \otimes\, (\text{det}_{b_m} \,W)^{\otimes n_m} \cong \text{det}_{b'} \,W\] for every $\b' \in\BB'$. The assignment $(W_i, w_a, f_\b) \mapsto  (W_i, w_a, f_{\b'})$ gives an equivariant isomorphism from $\RR(Q, l, \alpha)$ to $\RR(Q, l, \alpha)$, under which semistable points are sent to semistable points. This follows from the fact semistability depends only on the subgroup $R \subset \Wt(Q)$ and the factor $(W_i,w_a)$ of $W$. The factor $(W_i, w_a)$ is not altered by the proposed isomorphism; checking $\b(W)=0$ for basis elements $\b \in \BB$ is equivalent to checking $r(W)=0$ on every element $r \in R$. Hence the  equivariant isomorphism above defines an isomorphism of the resulting stacks $\MM_\theta(Q, l, \alpha)$.

This is not to say that the choice of basis is unimportant. It only becomes unimportant when we insist that the linear maps $f_\b: \kk \rightarrow \text{det}_b\,W$ are isomorphisms. Indeed, let $\b' = -\b$. Then given linear map $f_\b: \kk \rightarrow \text{det}_b\,W$ there exists a natural linear map $(f_\b)^\vee: \text{det}_{\b'}\,W \cong (\text{det}_b\,W)^\vee \rightarrow \kk^\vee \cong \kk$. If $f_\b$ is an isomorphism then we define $f_{\b'}:= (f_\b^\vee)^{-1}$, otherwise there is no natural definition for $f_{b'}$. In Section 5, we will allow the maps $f_b$ to be zero and will have to choose a basis carefully.
\end{rmk}

Given a quiver $Q$, a family of quiver representations $(\mathcal{W}_i, {w_a})$ over a base scheme $S$ and an element $\b = \sum_{i \in Q_0} \b_i \ee_i \in \ZZ^{Q_0}$ define
\begin{equation*}
\textup{det}_{\b}\, \mathcal{W}:= \bigotimes_{i\in Q_0}(\det \mathcal{W}_i)^{\otimes \b_i}.
\end{equation*}
Here we use the convention $\mathcal{W}_i ^{\otimes -1} := \mathcal{W}_i^{\vee}$. 

To justify the term `moduli stack' in Definition \ref{moduli} we must give a suitable notion of families over schemes for which the moduli stack is $\MM(Q, l, \alpha)$. One could define a family of refined representations over a scheme $S$ to be a refined representation of $(Q, l)$ in the category of locally free $\OO_S$-modules of $S$, that is, Definition \ref{fam} without the isomorphism of line bundles $\OO_S \rightarrow \text{det}_{\theta_\Delta}\,\mathcal{W}$. However, this would imply that $\Delta$ is a subgroup of the automorphism group of any given object. This gives stacks that are unsuitable for our applications as they don't admit closed immersions from Deligne-Mumford stacks. 

If $\alpha$ is primitive, that is the greatest common factor of its components is 1, we can alter the definition of a family to sidestep this issue. We do this by adding an extra nonzero parameter to $\RR(Q, l, \alpha)$ on which $\Delta$ acts with weight 1. This amounts to finding a character $\theta_\Delta \in \ZZ^{Q_0}$ for which $\langle \theta_\Delta, \Delta \rangle =\sum_i \theta_i \alpha_i= 1$, one may find such $\theta_\Delta$ precisely when $\alpha$ is primitive. For the rest of the section fix a primitive dimension vector $\alpha$ and pick $\theta_\Delta \in \ZZ^{Q_0}$ such that $\langle \theta_\Delta, \Delta \rangle = 1$.

\begin{defn}\label{fam}
A {\em flat family of refined representations of $(Q, l)$} over a connected scheme $S$ is a collection of rank $\alpha_i$ locally free sheaves $\mathcal{W}_i$ for $i \in Q_0$, together with a choice of morphisms $\mathcal{W}_{t(a)} \rightarrow \mathcal{W}_{h(a)}$ for $a \in Q_1$, isomorphisms of line bundles $\OO_S \rightarrow \text{det}_b\,\mathcal{W}$ for $\b \in \BB$ and an isomorphism of line bundles $\OO_S \rightarrow \text{det}_{\theta_\Delta}\,\mathcal{W}$.
\end{defn}

\begin{prop}
The stack $\MM_\theta(Q,l, \alpha)$ is the moduli stack of families of $\theta$-semistable refined representations  of $(Q, l)$.
\end{prop}

\begin{proof}
First we identify the nonzero elements of $\text{det}_{\theta_\Delta}\,W$ with $\kk^\times$, $\GL(\alpha)$ acts on $\text{det}_{\theta_\Delta}\,W$ by change of basis. Consider the stack quotient $[\RR(Q, l, \alpha)^{ss}_\theta \times \kk^\times / \GL(\alpha)]$. We claim that this represents the moduli problem defined by Definition \ref{fam}. An object in $[\RR(Q, l, \alpha)^{ss}_\theta \times \kk^\times / \GL(\alpha)](S)$ is a principal $\GL(\alpha)$-bundle $\mathcal{P}:= \bigoplus_{i\in Q_0} \mathcal{P}_i$ over $S$ with a $\GL(\alpha)$-equivariant morphism $\mathcal{P} \rightarrow \RR(Q, l, \alpha)^{ss}_\theta \times \kk^\times$. Define $\mathcal{W}_i$ to be the $W_i$-bundles corresponding to $\mathcal{P}_i$. Let $(U_j)_{j \in J}$ be an open cover of $S$ that trivializes $\mathcal{P}_i$. For $j \in J$ an equivariant morphism $U_j \times \GL(\alpha) \rightarrow \RR(Q, l, \alpha)^{ss}_\theta \times \kk^\times$ is determined by the image of the identity fibre and so is determined by a morphism $U_j \rightarrow \RR(Q, l, \alpha)^{ss}_\theta \times \kk^\times$. This morphism in turn defines a section of the vector bundle $U_j \times \Hom(W_{t(a)}, W_{h(a)})$ for every $a \in Q_1$, a nonzero section of $U_j \times \text{det}_b\,W$ for every $\b \in \BB$ and a nonzero section of $U_j \times \text{det}_{\theta_\Delta}\,W$. Since these sections come from a globally defined map $\mathcal{P} \rightarrow \RR(Q, l, \alpha)^{ss}_\theta \times \kk^\times$, they glue to give the required family over $S$. Similarly a family over $S$ defines an object of $[\RR(Q, l, \alpha)^{ss}_\theta \times \kk^\times / \GL(\alpha)](S)$. Morphisms of families correspond naturally to morphisms of objects of $[\RR(Q, l, \alpha)^{ss}_\theta \times \kk^\times / \GL(\alpha)](S)$.

The choice of $\theta_\Delta$ implies that $\Delta$ acts with weight one on the space of isomorphisms from $\kk$ to $\text{det}_{\theta_\Delta}\,W$. For every element $(W, t) \in \RR(Q, l, \alpha)^{ss}_\theta \times \kk^\times$ there exists a unique element of $\Delta$ that acts on $(W,t)$ to give $(W,1)$. The subgroup of $\GL(\alpha)$ that fixes the $\kk^\times$ component is isomorphic to $\PGL(\alpha)$, so we have a stack isomorphism 
\begin{equation*}
[\RR(Q, l, \alpha)^{ss}_\theta \times \kk^\times / \GL(\alpha)] \cong [\RR(Q, l, \alpha)^{ss}_\theta / \PGL(\alpha)] = \MM_\theta(Q, l, \alpha)
\end{equation*}
as required.
\end{proof}

The choice $\theta_\Delta$ might seem ad hoc at the moment, but a natural choice presents itself in our applications, see Remark \ref{conva}. 

\section{Quivers of sections}\label{Quiver}
In this section we introduce a generalization of the classical linear series construction to projective toric orbifolds, $\XX$. Starting with a collection of line bundles $\LL = (L_0, L_1,\ldots, L_r)$ on $\XX$ we produce a labelled quiver $(Q, \div)$ and give a rational map $\XX \dashrightarrow \MM_\theta(Q,\div, \alpha)$ with $\alpha =(1,\ldots,1)$. We then go on to study certain properties of this rational map. Throughout this section we assume that all our stacks $\XX$ are projective toric orbifolds; we will also fix the dimension vector to be $\alpha:= (1,\ldots,1)$ and drop it from the notation.

We now extend the definition of a quiver of sections, as defined by Craw-Smith in the beginning of Section 3 of \cite{Craw-Smith}, to toric orbifolds; we reproduce their definitions and adopt their conventions. Let $\LL: =  (L_0, L_1,\ldots, L_r)$ be a list of distinct line bundle on a stack $\XX$. A $T_\XX$-invariant section $s \in \Gamma(\XX, L_j \otimes L_i^\vee)$ is {\em irreducible} if it can not be written as a product of two nonzero sections $s' \in \Gamma(\XX, L_k \otimes L_i^\vee)$ and $s'' \in \Gamma(\XX, L_j \otimes L_k^\vee)$ for $L_k \in \LL$.  The {quiver of sections} associated to $\LL$ is a quiver $Q$ in which the vertices $Q_0 = \{0, \ldots, r\}$ correspond to the line bundles of $\LL$ and the arrows from $i$ to $j$ correspond to the set of irreducible $T_\XX$-invariant sections in $\Gamma(\XX, L_j \otimes L_i^\vee)$.

Since each arrow $a \in Q_1$ corresponds to a $T_\XX$-invariant section $s \in \Gamma(\XX, L_j \otimes L_i^\vee)$, we define a map $\div: Q_1 \rightarrow \ZZ^{\Sigma(1)}$ by sending $a$ to the corresponding divisor in $\ZZ^{\Sigma(1)}$. We call the labelled quiver $(Q, \div)$ the {\em labelled quiver of sections} of $\LL$.
\begin{conventions}\label{conv}
Let $(Q, \div)$ be the labelled quiver of sections corresponding to $\LL = (L_0, \ldots, L_r)$.
\begin{itemize}
\item[(a)] By definition, $(Q, \div)$ only depends on the line bundles $L_j \otimes L_i^\vee$ where $L_i, L_j \in \LL$. Consequently, for any line bundle $L'$ on $\XX$, we have $(Q, \div) = (Q', \div)$ where $(Q', \div)$ is a quiver of sections associated to $\LL' = (L_0 \otimes L', \ldots, L_r \otimes L').$ To eliminate this redundancy, we will assume that $L_0 = \OO_\XX$.
\item[(b)] We will assume that $\Gamma(\XX, L_i) \neq 0$ for $L_i \in \LL$. This implies $Q$ is connected and rooted at $0 \in Q_0$
\end{itemize}
\end{conventions}

\begin{rmk}\label{conva}
For labelled quivers of sections of line bundles, Convention \ref{conv} a) fixes $\theta_\Delta = (1, 0,\ldots, 0)$.
\end{rmk}

Keeping the notation of Section \ref{Quivers}, define $\pic: \Wt(Q) \rightarrow \Pic(\XX)$ by $\theta = \sum_{i \in Q_0} \theta_i \ee_i \mapsto \bigotimes_{i \in Q_0}L_i^{\otimes \theta_i}$ and let $\deg: \ZZ^{\Sigma(1)} \rightarrow \Pic(\XX)$ be the homomorphism in short exact sequence (\ref{Coxsq}). We then have the following commutative diagram
\begin{equation}\label{com1}\begin{split}
\xymatrix @C =1.5cm @R=1.3cm { \ZZ^{Q_1} \ar@{->>}[r]^-{\text{inc}} \ar[d]_{\text{div}} & \text{Wt}(Q) \ar[d]^{\text{pic}} \\
\ZZ^{\Sigma(1)} \ar@{->>}[r]^-{\text{deg}} & \text{Pic}(\XX).}
\end{split}
\end{equation}
The subgroup $R$ is by definition the image under inc of the kernel of div, so diagram (\ref{com1}) restricts to the following commutative diagram
\begin{equation}\label{com2}\begin{split}
\xymatrix @C =1.5cm @R=1.3cm { \,\,R\,\,\ar@{^{(}->}[r]^\iota \ar[d]_{0} & \text{Wt}(Q) \ar[d]^{\text{pic}} \\
\ZZ^{\Sigma(1)} \ar@{->>}[r]^-{\text{deg}} & \text{Pic}(\XX).}
\end{split}
\end{equation}

Define a $\Wt(Q)$-grading of the semigroup algebra $\kk[\NN^{Q_1} \oplus R]$ by assigning the monomial $y^{u}z^{v} \in \kk[\NN^{Q_1} \oplus R]$ degree $\inc(u) + \iota(v)$. This grading induces the change of basis action of $\PGL(\alpha)$ on $\RR(Q, \div) \cong \Spec(\kk[\NN^{Q_1} \oplus R])$. On the other hand, the map $\deg$ gives the $\Pic(\XX)$-grading of $\kk[\NN^{\Sigma(1)}]$ that arises from the short exact sequence (\ref{Coxsq}).

By construction $\div( \NN^{Q_1}) \subset \NN^{\Sigma(1)}$, so the map $\div \oplus 0: \NN^{Q_1} \oplus R \rightarrow \NN^{\Sigma(1)}$ induces a map of semigroup algebras $\Psi: \kk [\NN^{Q_1} \oplus R] \rightarrow \kk[ \NN^{\Sigma(1)}]$, which in turn defines a morphism $\Psi^*$ from $\AA^{\Sigma(1)}$ to $\RR(Q, \div)$. This morphism is equivariant with respect to the actions of the groups $\Hom(\Pic(\XX), \kk^\times)$ and $\PGL(\alpha) \cong \Hom(\Wt(Q),\kk^\times)$ on $\AA^{\Sigma(1)}$ and $\RR(Q, \div)$ because the diagrams (\ref{com1}) and (\ref{com2}) commute. Thus for any $\theta \in \Wt(Q)$, $\Psi$ induces a rational map 
\begin{equation*}
\psi_{\theta}: \XX \dashrightarrow \MM_\theta(Q, \div).
\end{equation*}
The rational map $\psi_\theta$ is a morphism of stacks \iff the inverse image under $\Psi^*$ of the $\theta$-unstable locus of $\RR(Q, \div)$ is contained in the Cox unstable locus $\VV(B_\XX)$, as defined in (\ref{Coxunstable}) (see Perroni's work on morphisms of toric stacks \cite{Perroni}).

We say a character $\chi_\theta \in \PGL(\alpha)^\vee$ is {\em generic} if every $\chi_\theta$-semistable point is $\chi_\theta$-stable.

\begin{defn} \label{bpf}
Take $\XX$ as above. A collection of line bundles $\LL = (\OO_\XX, L_1, \ldots, L_r)$ is {\em base-point free} if there exists generic $\chi_\theta \in \PGL(\alpha)^\vee$ for which $\psi_{\theta}: \XX \dashrightarrow \MM_\theta(Q, \div)$ is a morphism.
\end{defn}

\begin{rmk}
In the case where $\chi_\theta$ is not generic $\MM_\theta(Q,\div)$ is a not a Deligne-Mumford stack. Furthermore, it rarely has a coarse moduli space. We insist $\chi_\theta$ is generic to avoid such ambient stacks. 
\end{rmk}

For $L \in \Pic(\XX)$ let $(s_0, \ldots, s_n)$ be a basis of $\Gamma(\XX, L)$ and let $\varphi_{|L|}: \XX \dashrightarrow \PP(\Gamma(\XX, L)^\vee)$ be the rational map taking $x \in \XX$ to $(s_0(x), \ldots, s_n(x)) \in  \PP(\Gamma(\XX, L)^\vee)$. We will say $L$ is {\em base-point free} if $\varphi_{|L|}$ is a morphism. It can be shown using the universal property of a coarse moduli space, that every base-point free bundle on a stack can be pulled back from a base-point free line bundle on the coarse moduli space.

\begin{lemma}
Let $\XX$ be a projective toric orbifold. A nontrivial line bundle $L$ on $\XX$ is base-point free \iff the collection $\LL = (\OO_\XX, L)$ is base-point free.
\end{lemma}

\begin{proof}
The map $\pic$ takes the basis vector $\ee_1-\ee_0$ of $\Wt(Q)\cong \ZZ$ to $L$. This implies that $\ker(\pic)$ is trivial, otherwise $L^{\otimes n} \cong \OO_\XX$ for some $n>0$ contradicting projectivity of $\XX$.  Now $R$ is a subgroup of $\ker(\pic)$ and is therefore trivial. Therefore $\MM_\theta(Q,\div)$ is $\MM_\theta(Q)$. Since $Q$ is acyclic, the only chamber in $\Wt(Q)_\QQ \cong \QQ$ is $\QQ_{>0}$; take $\theta$ in this chamber, then $\MM_\theta(Q) \cong \PP(\Gamma(\XX, L))$ from which the claim follows.
\end{proof}

\begin{ex}\label{p112ex}
Let $\XX = \PP(1,1,2)$ and $\LL = (\OO_\XX, \OO(1), \OO(2))$. The labelled quiver of sections of $\LL$ is given by the quiver in Figure \ref{p112} and the labelling map $\div: \ZZ^{5} \rightarrow \ZZ^{\Sigma(1)} \cong \ZZ^{3}$ defined by the matrix
\begin{equation*}\begin{pmatrix}
1 & 0 & 1& 0&0 \\ 0&1&0&1&0 \\ 0&0&0&0&1
\end{pmatrix}.\end{equation*}
\begin{figure}[h!]
\centering
\begin{equation*}
\entrymodifiers={++[o][F-]}
\xymatrix @C=3pc{0 \ar@<0.8ex>[r]|{a_1} \ar@<-0.8ex>[r]|{a_2} \ar@/_-1.5pc/ [rr]^{a_5} & 1 \ar@<0.8ex>[r]|{a_3} \ar@<-0.8ex>[r]|{a_4} & 2}
\end{equation*}\caption{A quiver of section of $\PP(1,1,2)$.}\label{p112} 
\end{figure}

\noindent Since $\ee_{a_1} - \ee_{a_3}$ and $\ee_{a_2} - \ee_{a_4}$ generate $\ker(\div)$, the element $\ee_0 -2 \ee_1+\ee_2 \in \Wt(Q)$ generates $R$. The map $\Psi^*: \AA^3 \longrightarrow \AA^5 \times \kk^\times$ sends $(x_1,x_2,x_3)$ to $(x_1,x_2,x_1,x_2,x_3,1)$. Write $\AA^{5} \times \kk^\times \cong \Spec(\kk[y_{a_1},\ldots, y_{a_5}, z_f^{\pm}])$. For $\theta = -3\ee_0+2\ee_1+ \ee_2 \in \Wt(Q)$ the $\theta$-unstable locus is given by \[\VV(B_\theta):= \Big \langle y^uz^v \in \kk[\NN^{Q_1} \oplus \ZZ]\,\Big|\, \inc(u) + \iota(v) =\theta\Big\rangle.\] The set $(\inc \oplus\iota)^{-1}(\theta) \cap (\NN^{Q_1} \oplus \ZZ)$ contains the Laurent polynomials $$y_{a_1}^4z_f, y_{a_2}^4z_f, y_{a_3}^4z_f^{-3}, y_{a_4}^4z_f^{-3}, y_{a_5}^2z_f^{-1}.$$ Therefore $\VV(B_\theta) =\VV(y_{a_1},\ldots,y_{a_5})$. Hence, \[\MM_\theta(Q, \div) \cong [(\AA^{5} \times \kk^\times \setminus \{0\}\times \kk^\times) \,/\, (\kk^\times)^2] \cong \PP(1,1,1,1,2).\] Therefore $\Psi^*$ induces the map \[\psi_\theta: \PP(1,1,2) \xrightarrow{\hspace*{0.7cm}} \PP(1,1,1,1,2)\] that sends $(x_1,x_2,x_3)$ to $(x_1,x_2,x_1,x_2,x_3).$
\end{ex}

\begin{ex}\label{AHex1}
Let $\XX = \PP(1,1,2)$,  $\LL = (\OO_\XX, \OO(1), \OO(3))$. It can be shown that for $\theta = -2\ee_0 + \ee_1 + \ee_2$ the moduli stack $\MM_\theta(Q, \div)$ is isomorphic to $\PP(1,1,2,2,2,2)$ and \[\psi_\theta: \PP(1,1,2) \longrightarrow \PP(1,1,2,2,2,2)\]
is a morphism that sends $(x_1, x_2,x_3)$ to $(x_1,x_2,x_1^2,x_1x_2,x_2^2,x_3).$
\end{ex}

\begin{rmk}\label{AHex}
Example \ref{AHex1} recovers the Abramovich-Hassett \cite{AbramovichHassett} construction for $\XX = \PP(1,1,2)$ with polarizing line bundle $L= \OO(1)$ and natural numbers $n=0$ and $m=2$. More generally, given a polarizing line bundle $L$ on $\XX$, one may recover the Abramovich-Hassett construction when $n=0$ by applying our machinery to the collection  \[\LL = (\OO_\XX, L, L \otimes L^{\otimes 2}, \ldots, L^{\otimes m(m+1)/2})\]  if necessary, working with an `incomplete' quiver of sections. An incomplete quiver of sections is a quiver of sections where not all torus-invariant sections contribute to paths in the quiver, analogous to an incomplete linear series.
\end{rmk}

Let $L$ be a base-point free line bundle on a variety $X$. The linear series construction gives a morphism $\varphi_{|L|}: X \rightarrow \PP(\Gamma(X,L)^\vee)$, under which the pullback of the tautological line bundle on $\PP(\Gamma(X,L)^\vee)$ is $L$. The following proposition gives an analogous result.

\begin{prop}
Let $\LL = (\OO_\XX, L_1, \ldots, L_r)$ be base-point free. The pullback of the tautological bundles on $\MM(Q, \div)$ via $\psi_\theta$ is the collection $\LL$.
\end{prop}

\begin{proof}
The group $\GL(\alpha) \cong (\kk^\times)^{r+1}$ acts on $\text{det}_{\theta_\Delta}\,W \cong \Spec(\kk[y_{\theta_\Delta}])$ (cf. Remark \ref{conva}) by $(t_0, \ldots, t_r) \cdot y_{\theta_\Delta} = t_0 \cdot y_{\theta_\Delta}$. So the subgroup fixing nonzero $y_{\theta_\Delta}$ is given by $G_{\theta_\Delta}:= \{(t_0, \ldots, t_r) \in \GL(\alpha) \, |\, t_0 =1\}$. Restricting the projection map $\GL(\alpha) \twoheadrightarrow \PGL(\alpha)$ to $G_{\theta_\Delta}$ we get an isomorphism $\PGL(\alpha) \cong G_{\theta_\Delta}$.

The tautological line bundles of $\MM_\theta(Q, \div) \cong [(\RR(Q, l)^{ss}_\theta \times \kk^\times) / \GL(\alpha)]$ are given by the standard basis elements of $\ZZ^{Q_0} \cong \GL(\alpha)^\vee$. Under the isomorphism of stacks $[\RR(Q, l)^{ss}_\theta \times \kk^\times) / \GL(\alpha)] \cong [(\RR(Q, l)^{ss}_\theta / G_{\theta_\Delta}]$ the pullbacks of the tautological line bundles is given by the image of the basis elements of $\ZZ^{Q_0}$ under the map dual to the inclusion $G_{\theta_\Delta} \hookrightarrow \GL(\alpha)$; now under the isomorphism $\PGL(\alpha) \cong G_{\theta_\Delta}$ these are mapped to the elements $0, \ee_1-\ee_0, \ldots, \ee_r-\ee_0 \in \Wt(Q)$.

For $\eta \in \Wt(Q)$ the pullback of the associated line bundle of $[\RR(Q, l)^{ss}_\theta / \PGL(\alpha)]$ to $\XX$ is given by $\pic(\eta)$. Therefore the pullbacks of the tautological line bundles $0, \ee_1-\ee_0, \ldots, \ee_r-\ee_0 \in \Wt(Q)$ to $\XX$ are the line bundles $\OO_\XX, L_1, \ldots, L_r$ as required.
\end{proof}

We now work towards a condition on a collection $\LL$ that guarantees it is base-point free. 

\begin{lemma}\label{tor}
Let $\XX$ be a projective toric orbifold and $L_1, \ldots,L_n \in \Pic(\XX)$ be such that every $L_i$, for $1\leq i\leq n$, is base-point free. Given a section $s$ of $L:=L_1 \otimes \cdots \otimes L_n$, there exists $m \in \NN$ such that $s^m$ is in the image of the multiplication map \[\mu_m: \Gamma(\XX,L_1)^{\otimes m} \otimes_\kk \cdots \otimes_\kk \Gamma(\XX,L_n)^{\otimes m} \rightarrow \Gamma(\XX, L^{\otimes m}).\]
\end{lemma}

\begin{proof}
Let $\kk[x_0,\ldots,x_n]$ be the Cox ring of $\XX$ and $\mu: \Gamma(\XX, L_1) \otimes \cdots \otimes \Gamma(\XX,L_n) \rightarrow \Gamma(\XX,L)$ be the multiplication map. The line bundles $L_i$ are base-point free, they are therefore pullbacks of line bundles on the underlying coarse moduli space, which is a toric variety, and so correspond to polytopes $P_{L_i}$. The polytope $P_L$ corresponding to $L$ is given by $P_{L_1} + \ldots + P_{L_n}$ (see page 69 of Fulton \cite{Fulton}). While the lattice points of $P_L$ are not, in general, a sum of lattice points of the polytopes $P_{L_i}$, the vertices of $P_L$ are given by sums of vertices of $P_{L_i}$ (see page 11 of Sturmfels \cite{Sturmfels}). Therefore the sections corresponding to the vertices of $P_L$ lie in $\text{im}(\mu)$. Since the vanishing locus of the sections corresponding to the vertices of $P_L$ is equal to that of the sections of $L$, we have that the vanishing locus of the sections in $\text{im}(\mu)$ is equal to that of the sections of $L$. Now let $s \in \Gamma(\XX, L)$, then by Hilbert's Nullstellensatz there exists a natural number $m \in \ZZ$ such that $s^m$ is in the ideal generated by $\text{im}(\mu)$ as required.
\end{proof}

For a collection of line bundles $\LL = (\OO_\XX, L_1, \ldots, L_r)$ on $\XX$, define 
\begin{equation*}
\LL_{\text{bpf}}:= \{L_i^\vee \otimes L_j\, |\, L_i, L_j \in \LL \text{ and } L_i^\vee \otimes L_j \text{ is base-point free}\}
\end{equation*}
and let $\pic_\QQ: \Wt(Q)_\QQ \rightarrow \Pic(\XX)_\QQ$ be $\pic \otimes \text{id}$.

\begin{lemma}\label{lem}
If $\rank(\ZZ \LL) = \rank(\ZZ \LL_{\textup{bpf}})$ then $\ker(\pic_\QQ) \subset R_\QQ$.
\end{lemma}

\begin{proof}
In this proof we use additive notion for the binary operation on the Picard group to avoid confusion with $- \otimes \QQ$. 

Let $\omega = \sum_{i\in Q_0} \omega_i \otimes q_i \in \ker(\pic_\QQ)$ and pick $n \in \NN$ sufficiently large so that $n\, \omega = \sum_{i \in Q_0} n_i\, \omega_i \otimes 1$ for $n_i \in \ZZ$ and set $\lambda := \sum_{i \in Q_0} n_i \,\omega_i$. Since $R_\QQ$ is a vector subspace of $\Wt(Q)_\QQ$ we have $\omega \in R_\QQ$ \iff $n  \,\omega \in R_\QQ$ for any $n \in \QQ \setminus \{0\}$. Therefore it suffices to show $n \,\omega = \lambda \otimes 1 \in R_\QQ$, that is there exists an element $\tau \in \ZZ^{Q_1} \otimes 1$ such that $\inc_\QQ(\tau) =\lambda \otimes 1$ and $\div_\QQ(\tau) = 0$.

Take the basis $E:= \{\ee_1 - \ee_0, \ldots, \ee_r-\ee_0\}$ of $\Wt(Q)$ and write $\lambda$ as a difference of positive and negative parts, that is, write $\lambda = \lambda_+ - \lambda_-$ for $\lambda_+, \lambda_- \in \NN E$ without cancellation. Let $L_\pm = \pic(\lambda_\pm)$. The fact that $\lambda \otimes 1 \in \ker(\pic_\QQ)$ implies $L_+ \otimes 1 = L_- \otimes 1$ and the rank assumption gives us
\begin{equation}\label{q}
L_+ \otimes 1 = L_- \otimes 1 = \sum L_{b_i} \otimes q_i \quad  \text{for } L_{b_i} \in \LL_\text{bpf} \text{ and } q_i \in \QQ.
\end{equation}
We may take $n$ big enough to ensure that each $q_i \in \ZZ$. Rearrange equations (\ref{q}) to get
\begin{equation}\label{+}
(L_+ \otimes 1) + \bigg(\sum_{q_i <0} -q_i L_{b_i} \otimes 1\bigg) = \sum_{q_i>0} q_i L_{b_i} \otimes 1
\end{equation}
\begin{equation}\label{-}
(L_- \otimes 1) + \bigg(\sum_{q_i <0} -q_i L_{b_i} \otimes 1\bigg) = \sum_{q_i>0} q_i L_{b_i} \otimes 1.
\end{equation}
Fix a section of each of the following line bundles: $L_+$, $L_-$ and $L_{b_i}$ for which $q_i<0$ (in turn fixing a section of $\sum_{q_i < 0} -q_i L_{b_i}$). Using equations (\ref{+}) and (\ref{-}), this fixes sections $s_\pm$ of $\sum_{q_i>0} q_i L_{b_i} + L_{t\pm}$ for some torsion line bundles $L_{t\pm}$. Without loss of generality we assume $L_{t\pm} =0$, otherwise multiply $n$ in the beginning of the proof by the orders of $L_{t\pm}$.

Since the incidence map is onto $\Wt(Q)$ there exists elements $\tau_1 \in \ZZ^{Q_1}$ and $\tau_2 \in \ZZ^{Q_1}$ such that $\inc(\tau_1) = \lambda_+$ and $\inc(\tau_2)= \lambda_-$. By Lemma \ref{tor} there exists $m_\pm \in \NN$ such that $s_\pm^{m_\pm}$ are a product of sections of the line bundles $L_{b_i}$ for which $q_i>0$. By definition of the quiver of sections $Q$, every section of a line bundle in $\LL_\text{bpf}$ gives rise to a path in the quiver, so there exists $\tau_\pm \in \ZZ^{Q_1}$ such that $\div(\tau_\pm)= s_\pm^{m_\pm}$. Define 
\begin{equation*}
\tau:= (\tau_1 \otimes 1) - \bigg(\tau_+ \otimes \frac{1}{m_+}\bigg) - (\tau_2 \otimes 1) + \bigg(\tau_- \otimes \frac{1}{m_-}\bigg).
\end{equation*}
We have that $\inc_\QQ(\tau) = \lambda \otimes 1$ because $\inc_\QQ(\tau_+ \otimes \frac{1}{m_+}) = \inc_\QQ(\tau_- \otimes \frac{1}{m_-})$. We also have that $(\tau_1 \otimes 1) - (\tau_+ \otimes \frac{1}{m_+})$ and $(\tau_2 \otimes 1) - (\tau_- \otimes \frac{1}{m_-})$ map via div to the section of the line bundle $\sum_{q_i < 0} -q_i L_{b_i}$ fixed above, and hence $\div_\QQ(\tau) = 0$, as required.
\end{proof}

The following example highlights that tensoring with $\QQ$ in the statement of Lemma \ref{lem} is necessary.

\begin{ex}
Consider the $(\ZZ \oplus \ZZ/2\ZZ \oplus \ZZ/2\ZZ)$-grading of $\kk[x_1,x_2,x_3,x_4]$ given by:
\begin{equation*}
\deg(x_1) = (1, 0, 0);\,\, \deg(x_2)=(1,1,0);\,\, \deg(x_3)= (1,0,1);\,\, \deg(x_4)=(1,1,1)
\end{equation*}
and let $(\kk^\times \times\ZZ/2\ZZ\times \ZZ/2\ZZ) \curvearrowright \AA^4$ be the corresponding action. Take \[\XX = [(\AA^4\setminus \{0\}) / (\kk^\times \times\ZZ/2\ZZ\times \ZZ/2\ZZ)]\] and 
\begin{equation*}
\LL = (\OO, \OO(1,0,0), \OO(1,1,0), \OO(1,0,1), \OO(1,1,1),\OO(2,0,0)).
\end{equation*}
We have that $\OO(2,0,0)$ is base-point free and so is an element of   $\LL_\text{bpf}$ therefore $\rank(\ZZ \LL) = \rank(\ZZ \LL_{bpf})$. However it can be shown that $\ker(\pic) \nsubseteq R$.
\end{ex}

\begin{thm}\label{prop}
If $\rank(\ZZ \LL) = \rank(\ZZ \LL_\textup{bpf})$ then $\LL$ is base-point free.
\end{thm}

\begin{proof}
Let $\theta \in \Wt(Q)$. The associated character $\chi_\theta \in \PGL(\alpha)^\vee$ gives a morphism of stacks $\XX \rightarrow \MM_\theta (Q, \div)$ \iff the inverse image of the $\theta$-unstable points of $\RR(Q, \div)$ is contained in $\VV(B_\XX)$. After picking a higher multiple if necessary, we may assume that the $\theta$-unstable locus in $\RR(Q, \div)$ is precisely the vanishing locus of the monomial ideal 
\begin{equation*}
B_\theta:= \Big\langle\, y^uz^v \in \kk[\NN^{Q_1} \oplus R] \, \Big | \, \inc(u) + \iota(v) = \theta\,\Big\rangle 
\end{equation*}
and its inverse image $\psi^{-1} ( \VV(B_\theta)) \subset \AA^{\Sigma(1)}$ is the vanishing locus of the monomial ideal
\begin{equation*}
\div B_\theta:= \Big\langle\, x^{\div((u,v))} \in \kk[\NN^{\Sigma(1)}] \, \Big | \, \inc(u) + \iota(v) = \theta\,\Big\rangle.
\end{equation*}
Now let $L:= \pic(\theta)$ and define 
\begin{equation*}
B_L:= \Big\langle\, x^\nu \in \kk[\NN^{\Sigma(1)}] \, \Big| \, \deg(\nu)=L\,\Big\rangle.
\end{equation*}
Given $\chi_\theta$ for which $L= \pic(\theta) \in \LL_{\textup{bpf}}$, pick $m$ big enough such that the $m\theta$-unstable locus is cut out by $B_{m\theta}$. The line bundle $L \in \LL_{\textup{bpf}}$ so for every $\nu \in \NN^{\Sigma(1)}$ for which $\deg(\nu)= L$ there exists an element $\rho_s \in \ZZ^{Q_1}$ such that $\div(\rho_s,0)=\nu$, so $B_L = \div B_\theta$. The ideal $\div B_{m\theta}$ is contained in $B_{L^{\otimes m}}$ and the vanishing locus of $\div B_{m\theta}$ is contained in that of $\div B_\theta$. Therefore \[\VV(B_{L^{\otimes m}}) \subset \VV(\div B_{m \theta}) \subset \VV(\div B_\theta) = \VV(B_L).\] Since $\pic(m\theta)=L^{\otimes m}$ and $L$ is base-point free, Lemma \ref{tor} implies $\VV(B_L) = \VV(B_{L^{\otimes m}})$, therefore $\VV(\div B_{m \theta}) = \VV(\div B_\theta)$. The line bundle $L$ base-point free so $\VV(B_L) \subset \VV(B_\XX)$ and hence $\psi^{-1} (\VV(B_{m\theta}))  \subset \VV(B_\XX)$. Hence,  the rational map $\psi_{\theta}: \XX \dashrightarrow \MM_\theta(Q, \div)$ is in fact a morphism of stacks.

It remains to show we may pick a generic character for which $\psi_{\theta}: \XX \rightarrow \MM_\theta(Q, \div)$ is a morphism.
Let $S := \{\ee_j - \ee_i \in \Wt(Q) \,| \, \pic(\ee_j - \ee_i)\text{ is base-point free} \}$ and $\sigma \subset \Wt(Q)_\QQ$ be the cone generated by elements of $S$ and $R$. Since the generators of the cone map under pic to base-point free line bundles, we have that any $\theta$ in $\sigma$ gives a morphism $\psi_{\theta}: \XX \rightarrow \MM_\theta(Q, \div).$ We claim that $\sigma \subset \Wt(Q)_\QQ$ is top dimensional. The vector space $\Wt(Q)_\QQ$ is isomorphic to  $(\text{ker}(\pic_\QQ) )\oplus(\text{im}(\pic_\QQ))$. The image of $\pic_\QQ$ is generated by $\LL$, the rank assumption then implies that the elements of $\LL_\text{bpf}$ also generate $\text{im}(\pic_\QQ)$. We have $\pic(S) = \LL_\text{bpf}$. In addition Lemma \ref{lem} give us that $\ker(\pic _ \QQ) \subset R _ \QQ$, therefore elements of $\sigma$ span $\Wt(Q)_\QQ$ which proves the claim. So one may pick a generic $\theta$ in the interior of $\sigma$ that gives a well defined morphism, hence $\LL$ is base-point free.
\end{proof} 

The next example gives a base-point free collection $\LL$ that does not satisfy the rank condition $\rank(\ZZ \LL) = \rank(\ZZ \LL_\text{bpf})$.

\begin{ex}
Let $\XX = \PP(1,2,3)$ and $\LL = (\OO_\XX, \OO(1), \OO(2), \OO(3))$. Note that $\LL_{\text{bpf}}$ is empty. It can be shown that given $\theta \in \NN(\ee_1-\ee_0) \oplus \NN(\ee_2-\ee_0) \oplus \NN(\ee_3-\ee_0)$, $\MM_\theta(Q, \div) \cong \PP(1,1,1,2,2,3)$ and we have a morphism
\[\psi_\theta: \PP(1,2,3) \xrightarrow{\hspace*{0.7cm}} \PP(1,1,1,2,2,3)\] that sends $(x_1, x_2,x_3)$ to $(x_1,x_1,x_1,x_2,x_2,x_3)$.
\end{ex}

Given a base-point free collection of line bundles $\LL$, the next proposition explicitly describes the image of $\psi_\theta$. Let $I_\LL \subset \kk[\NN^{Q_1} \oplus R]$ be the ideal given by the following
\begin{equation}\label{IL}
I_\LL := \Big\langle y^{u_1}z^{v_1} - y^{u_2}z^{v_2} \, | \, \div(u_1-u_2) = 0, \inc(u_1-u_2)+\iota(v_1-v_2)=0 \Big\rangle.
\end{equation}

\begin{prop}\label{image}
Let $\LL$ be a base-point free collection of line bundles on $\XX$ and $\theta \in \Wt(Q)$ be such that $\psi_\theta$ is a morphism. Then the image of $\psi_\theta$ is given by $[(\VV(I_\LL) \setminus \VV(B_\theta)) / \PGL(\alpha)] \subset \MM_\theta(Q,\div)$.
\end{prop}

\begin{proof}
The image of the map from $\AA^{\Sigma(1)}$ to $\AA^{Q_1} \times (\kk^\times)^R$ induced by the semigroup homomorphism $\div \oplus 0: \NN^{Q_1} \oplus R \rightarrow \NN^{\Sigma(1)}$ is given by the vanishing locus of the toric ideal
\begin{equation}
I:= \Big\langle y^{u_1}z^{v_1} - y^{u_2}z^{v_2} \in \kk[\NN^{Q_1} \oplus R] \, \Big| \, \div(u_1-u_2) = 0 \Big\rangle.
\end{equation}
For any element $y^{u_1}z^{v_1} - y^{u_2}z^{v_2} \in \kk[\NN^{Q_1} \oplus R]$ , its $\Wt(Q)$-grade is defined to be $\inc(u_1-u_2)+\iota(v_1-v_2)$. Therefore, $I_\LL \subset \kk[\NN^{Q_1} \oplus R]$ is defined to give the $\Wt(Q)$-homogenous part of $I$. We conclude that the image $\psi_\theta$ is given by $$[(\VV(I_\LL) \setminus \VV(B_\theta)) / \PGL(\alpha)].$$
\end{proof}

\begin{rmk}
Let $\theta_1, \theta_2 \in \Wt(Q)$ be generic and such that the maps $\psi_{\theta_1}, \psi_{\theta_2}$ are morphisms. The fact that $\psi_{\theta_1}, \psi_{\theta_2}$ are morphisms implies that their images do not intersect the unstable loci $\VV(B_{\theta_1})$ and $\VV(B_{\theta_2})$ and so are independent of the unstable loci. Since the only difference between the two morphisms is the unstable loci of the target this implies that the image of $\psi_{\theta_1}$ is isomorphic to that of $\psi_{\theta_2}$.
\end{rmk}

Now we investigate representability of the morphism $\psi_\theta$. Let $\pi: \XX \rightarrow X$ be the map to the \cms $X$ of $\XX$. We recall a useful definition from Nironi \cite{Nironi}.
\begin{defn}[Def. 2.2, \cite{Nironi}]
A locally free sheaf $\mathcal{V}$ on $\XX$ is {\em $\pi$-ample} if for every geometric point of $\XX$ the representation of the stabilizer group at that point is faithful.
\end{defn}

\begin{thm}\label{rep}
Let $\LL= (\OO_\XX, L_1, \ldots, L_r)$ be a \bsf collection of line bundles. Then $\psi_{\theta}$ is representable \iff $\bigoplus_{j=1}^rL_j$ is $\pi$-ample.
\end{thm}

\begin{proof}
By Lemma 2.3.9 of \cite{AbramovichHassett}, $\psi_{\theta}$ is representable \iff  the map $g: \Aut(x) \rightarrow \Aut(\psi_{\theta}(x))$ is injective for every geometric point $x \in \XX$.  The map $g$ fits into the following commutative diagram
\begin{equation}\label{com4}\begin{split}
\xymatrix@C =1.3cm @R=1.3cm{\text{Aut}(x) \ar[r]^-{g} \ar@{^{(}->}[d] & \text{Aut}(\psi_\theta(x)) \ar@{^{(}->}[d] \\ \Hom(\Pic(\XX), \kk^\times) \ar[r]^{\pic^\vee} & \Hom(\Wt(Q), \kk^\times).}
\end{split}\end{equation}
Here $\pic^\vee$ denotes the map given by applying the functor $\Hom(-,\kk^\times)$ to $\pic$. We claim that the representation of $\Aut(x)$ given by $\bigoplus_{j=1}^rL_j$ is the composite \begin{equation}\label{blah} \Aut(x) \hookrightarrow \Hom(\Pic(\XX), \kk^\times) \xrightarrow{ \pic^\vee }  \Hom(\Wt(Q), \kk^\times).\end{equation} Indeed, take the basis $\{\ee_i - \ee_0 \in \Wt(Q) \,|\, i= 1,\ldots r\}$ of $\Wt(Q)$ giving an isomorphism $ \Hom(\Wt(Q), \kk^\times) \cong (\kk^\times)^r$. Evaluating at $\ee_i-\ee_0 \in \Wt(Q)$ gives a map $\Hom(\Wt(Q), \kk^\times) \rightarrow \kk^\times$. By definition of the $\Hom$-functor the composite \[\Hom(\Pic(\XX), \kk^\times) \xrightarrow{ \pic^\vee }  \Hom(\Wt(Q), \kk^\times) \rightarrow \kk^\times\] is given by evaluating at $\pic(\ee_i - \ee_0) = L_i$ and is therefore the representation induced by the line bundle $L_i$. This proves the claim. So we have that the composite (\ref{blah}) is injective for every $x \in \XX$ precisely when $\bigoplus_{j=1}^rL_j$ is $\pi$-ample. Commutativity of (\ref{com4}) gives that (\ref{blah}) is injective \iff $g$ is injective, as required.
\end{proof}

For $\XX$ a toric orbifold, let $\pi: \XX \rightarrow X$ be the map to the coarse moduli space. The group homomorphism given by the pullback functor, $\pi^* \colon \Pic(X) \rightarrow \Pic(\XX)$, identifies $\Pic(X)$ with a subgroup of $\Pic(\XX)$. We will abuse notation and use $\Pic(X) \subset \Pic(\XX)$ to denote this subgroup.

\begin{coro}
Let $\LL$ be a base-point free collection of line bundles on $\XX$. If $\LL$ generates $\Pic(\XX) / \Pic(X)$ then $\psi_\theta$ is representable.
\end{coro}

\begin{proof}
This follows from the fact that elements of $\Pic(X)$ give trivial representations of $\Aut(x)$ for every geometric point $x\in \XX$.
\end{proof}

The following example shows that for a given base-point free collection $\LL$, representability of $\psi_\theta$ is weaker than $\LL$ generating $\Pic(\XX) / \Pic(X)$.
\begin{ex}
Take $N= \ZZ$ and let $\Sigma$ be the fan associated to the toric variety $\PP^1$ with rays $\rho_+ := \QQ_{\geq 0}$ and $\rho_-:= \QQ_{\leq 0}$. Let $\beta: \ZZ^{\Sigma(1)} \rightarrow \ZZ$ take $\ee_{\rho_{\pm}}$ to $\pm2$. For $\mathbf{\Sigma} = (N, \Sigma, \beta)$ take $\XX = \XX_{\mathbf{\Sigma}}$. First note that $\Pic(\XX) \cong \ZZ \oplus \ZZ/2\ZZ$ and let $\LL= (\OO_\XX, \OO(2,1), \OO(4,0))$. For $\theta = -3\ee_0 + 2\ee_1+\ee_2 \in \Wt(Q)$ the collection $\LL$ gives a representable morphism $\psi_\theta: \XX \rightarrow \PP(1,1,2,2)$ taking $(x,y) \in \AA^{\Sigma(1)}$ to $(xy, xy, x^4, y^4)$.
\end{ex}

\begin{ex}
Every example in this section defines a representable morphism. For an example that does not, take $\XX = \PP(1,1,2)$ and $\LL=(\OO, \OO(2))$.
\end{ex}

We conclude the section by comparing our construction to that of Craw-Smith \cite{Craw-Smith} in the case where $\XX=X$ is a toric variety. Let  $\LL = (\OO_X, L_1, \ldots, L_r)$ be a collection of base-point free line bundles on a toric variety $X$ and $\vartheta =-r \ee_0+ \ee_1+\cdots+\ee_r$, then Craw and Smith use the commutative diagram (\ref{com1}) to produce a morphism $\varphi_{|\LL|}: X \rightarrow \MM_\vartheta(Q)$.

\begin{prop}
Let $\XX = X$ be a toric variety and $\LL$ be a collection of \bsf line bundles on $X$. Then the image of $\psi_{\vartheta}$ is isomorphic to that of $\varphi_{|\LL|}$.
\end{prop}

\begin{proof}
We have the following commutative diagram:
\begin{equation}\begin{split}
\label{comdia}
\xymatrix @R=1cm @C=2.5cm {\ZZ^{Q_1} \ar[r]^-{\text{inc}}\ar[d]_-{(\text{id},0)} & \text{Wt}(Q) \ar[d]^{\text{id}} \\
\ZZ^{Q_1} \oplus R  \ar[r]^{\text{inc} \oplus \iota} \ar[d]_{\text{div}\oplus 0} & \text{Wt}(Q) \ar[d]^{\pic} \\
\ZZ^{\Sigma(1)} \ar[r]^-{\text{deg}} & \text{Pic}(X).}\end{split}
\end{equation}
The maps of semigroups $\NN^{Q_1} \xrightarrow{(\text{id}, 0)} \NN^{Q_1} \oplus R \xrightarrow{\div\oplus 0} \NN^{\Sigma(1)}$ give maps of semigroup algebras $\kk[\NN^{Q_1}] \rightarrow \kk[\NN^{Q_1} \oplus R] \rightarrow \kk[\NN^{\Sigma(1)}]$. After applying the functor Spec these give morphisms
\begin{align*}
\AA^{\Sigma(1)} \xrightarrow{\hspace{0.2cm}\Psi^*\hspace{0.2cm}} \AA^{Q_1} \times (\kk^\times)^{P} \xrightarrow{\hspace{0.2cm}\pi\hspace{0.2cm}} \AA^{Q_1}.
\end{align*}
The morphism $\Psi^*$ is induced by the semigroup map $\NN^{Q_1} \oplus R \xrightarrow{\div\oplus 0} \NN^{\Sigma(1)}$, therefore its image lies in the subvariety $\AA^{Q_1} \times (1,\ldots,1) \subset \AA^{Q_1} \times (\kk^\times)^P$. On the other hand, the morphism $\pi$ is just the projection to the first factor, therefore the image of $\Psi^*$ is isomorphic to the image of $\pi \circ \Psi^*$.
The commutativity of diagram (\ref{comdia}) implies that $\Psi^*$ is equivariant with respect to the actions of the groups $\Hom(\Pic(X), \kk^\times)$ and $\Hom(\Wt(Q), \kk^\times)$ induced by $\deg$ and $\inc$, similarly $\pi$ is equivariant. So the maps $\Psi^*$ and $\pi$ give rise to rational maps:
\begin{equation*}
\xymatrix{X \ar@{-->}[r]^-{\psi_\vartheta} & \MM_\vartheta(Q,\div) \ar@{-->}[r]^-\pi &\MM_\vartheta(Q)}.
\end{equation*}
The composite $\pi \circ \psi_\theta$ is equal to the rational map $\varphi_{|L|}$ and since $\LL$ is a collection base-point free line bundles Corollary 4.2 of \cite{Craw-Smith} implies it is a morphism. By virtue of Proposition \ref{prop} we have that $\psi_\vartheta$ is a morphism. Now, by definition, $\psi_\vartheta$ and $\varphi_{|L|}$ descend from $\Psi^*$ and $\pi \circ \Psi^*$ respectively and since $\Psi^*$ and  $\pi \circ \Psi^*$ have isomorphic images the images of $\psi_\vartheta$ and $\varphi_{|L|}$ are isomorphic.
\end{proof}

\section{Application to the McKay correspondence}\label{McKay}
In this section, we apply the construction from the previous section to toric quotient singularities. For $G \subset \GL(n, \kk)$ a finite abelian group and $(Q, \div)$ the labelled McKay quiver, we construct a closed immersion $[\AA^n / G] \hookrightarrow \MM_\theta(Q, \div)$ for any $\theta \in \Wt(Q)$. In the case where $G \subset \SL(n, \kk)$ and $n \leq 3$, we proceed to alter the construction of $\MM_\theta(Q, \div)$ and yield a GIT problem for which one generic stability condition gives $[\AA^n / G]$ and another gives $G$-Hilb$(\AA^n)$. We will assume $\alpha = (1,\ldots,1)$ throughout this section and drop it from the notation.

Take $n \in \NN$ and $G$ a finite abelian subgroup of $\GL(n, \kk)$ with no quasireflections. We may assume that $G$ is contained in the subgroup $(\kk^\times)^n$ of diagonal matrices with nonzero entries in $\GL(n, \kk)$. Line bundles on $[\AA^n/G]$ are given by $G$-equivariant line bundles on $\AA^n$, which in turn are determined by $G$-equivariant isomorphisms $\OO_{\AA^n \times G} \rightarrow \OO_{\AA^n \times G}$. From this it follows that the Picard group  of $[\AA^n/G]$ is naturally isomorphic to the group characters $G^\vee$. With these preparations, take \[\LL = (\OO_{\AA^n} \otimes \rho \,|\, \rho \in G^\vee).\] Then the labelled quiver of sections $(Q,\div)$ of $\LL$ coincides with the McKay quiver, see the beginning of Section 4.1 of \cite{CQV}. 

From now on we will use the isomorphism $\Pic([\AA^n / G]) \cong G^\vee$ tacitly. In much the same way as we have commutative diagrams (\ref{com1}) and (\ref{com2}) we have
\begin{equation} \label{com5} \begin{split}
\xymatrix @C =1.5cm @R=1.3cm { \ZZ^{Q_1} \ar@{->>}[r]^-{\text{inc}} \ar[d]_{\text{div}} & \text{Wt}(Q) \ar[d]^{\text{pic}} & \,\,R\,\, \ar@{^{(}->}[r]^{\iota} \ar[d]_{0} & \Wt(Q) \ar[d]^{\pic} \\
\ZZ^n \ar[r]^-{\text{deg}} & G^\vee & \ZZ^{Q_1} \ar[r]^-{\deg} & G^\vee.}
\end{split}\end{equation} 
As in Section \ref{Quiver}, the semigroup morphism $\div \oplus 0: \NN^{Q_1} \oplus R \rightarrow \NN^n$ gives a morphism $\Psi^*:\AA^n \rightarrow \AA^{Q_1} \times (\kk^\times)^R$. The commutativity of (\ref{com5}) gives that $\Psi^*$ is equivariant with respect to the actions of $G$ on $\AA^n$ and $\Hom(\Wt(Q), \kk^\times)$ on $\AA^{Q_1} \times (\kk^\times)^R$. Given a $\theta \in \Wt(Q)$ this gives a rational map \[\psi_\theta: [\AA^n/G] \dashrightarrow \MM_\theta(Q, \div).\]

From now on we identify the lattice $\Wt(Q)$ with the lattice $\{\theta \in \ZZ^{Q_0} \, |\, \theta_0 =0\}$ whose basis is $\{\ee_{\rho} \,|\, \rho \in G^\vee \setminus \{0\}\}$. 

\begin{prop}\label{clsdimm}
For any $\chi_\theta \in \PGL(\alpha)^\vee$, \[\psi_\theta: [\AA^n/G] \dashrightarrow \MM_\theta(Q, \div)\] is a closed immersion.
\end{prop}
\begin{proof}
We begin by studying the $\theta$-semistable points. By definition of the quiver of sections, a path $p$ from $\rho_0$ to $\rho$ corresponds to a section $s \in \Hom(\rho_0, \rho)$. Then there exists a path $p'$ from $\rho'$ to $\rho \otimes \rho'$ given rise to by the same section $s \in \Hom(\rho', \rho' \otimes \rho)$, that is $\div(p) = \div(p')$. Note that $\inc(p) = \ee_\rho $ and $\inc(p')= \ee_{\rho' \otimes \rho} - \ee_{\rho'}$. Since $\div(p') - \div(p) =0$ we have that $ -\ee_\rho  - \ee_{\rho'}   + \ee_{\rho\otimes \rho'} \in R$. Given that $\ker(\pic)$ is generated by elements of the form $-\ee_\rho- \ee_{\rho'} + \ee_{\rho \otimes \rho'}$ this shows $\ker(\pic) \subset R$. The commutativity of the diagrams (\ref{com5}) gives $R \subset \ker(\pic)$ and therefore $R = \ker(\pic)$. 

The image of $\pic$ is a torsion $\ZZ$-module, so $R = \ker(\pic)$ $\QQ$-spans $\Wt(Q)$ and hence any basis $\BB$ of $R$ $\QQ$-spans $\Wt(Q)$. Now let $W$ be a refined representation and let $\theta \in \Wt(Q)$. We claim that $W$ is $\theta$-semistable. Indeed, let $W_\bullet$ be $\kk Q$-module filtration satisfying the conditions in Definition \ref{stab} and write $\theta$ as a $\QQ$-linear combination of $\b \in \BB$. Then since $\b(W_\bullet)=0$, we have $\theta(W_\bullet)=0$. This in particular implies that the $\theta$-unstable locus is empty. It follows at once that the rational map $\psi_\theta$ is a morphism \[\psi_\theta: [\AA^n/ G] \longrightarrow \MM_\theta(Q, \div)\] for any $\theta \in \Wt(Q)$. Let $I_\LL$ be the $\kk[\NN^{Q_1} \oplus R]$ ideal defined in (\ref{IL}). Again, after noting the $\theta$-unstable locus is empty, an argument similar to that of Proposition \ref{image} gives that the image of $\psi_\theta$ is $[\VV(I_\LL) / \PGL(\alpha)]$. It remains to show that $[\AA^n /G]$ is isomorphic to $[\VV(I_\LL) / \PGL(\alpha)]$. Consider $\kk[\NN^{Q_1} \oplus R] / I_\LL$ and multiply the generators $y^{u_1}z^{v_1} - y^{u_2}z^{v_2}$ of $I_\LL$ by the units $z^{-v_1}$ to get an alternative set of generators $y^{u_1}-y^{u_2}z^{v_2-v_1}$. Then $I_\LL$ is given by \[ \Big \langle y^{u_1}-y^{u_2}z^v \in \kk[\NN^{Q_1} \oplus R] \,\Big|\, \div(u_1 - u_2) =0, \inc(u_1 -u_2)- \iota(v) =0 \Big\rangle. \] 
Pick $a_1, \ldots, a_n \in Q_1$ such that $\div(a_i)$ is the $i$th basis element of $\ZZ^n$. Since every arrow in the McKay quiver is labelled by a basis element of $\ZZ^n$ the kernel of $\div$ is generated by differences $\ee_{a_i'} -\ee_{a_i}$ with $\div(a_i') = \div(a_i)$. By definition of $R$, for every generator $\ee_{a_i'} -\ee_{a_i}$ of $\ker(\div)$ there exists $v' \in R$ such that $\inc(\ee_{a_i} -\ee_{a_i'}) = v'$. Therefore $I_\LL$ is generated by elements of the form $y_{a_i'} - y_{a_i}z^{v'}$. This implies that for every $a \in Q_1$ not in the list $a_1, \ldots, a_n$, the monomial $y_a$ is equivalent in the quotient $\kk[\NN^{Q_1} \oplus R] / I_\LL$ to a product of elements in $\kk[y_{a_1}, \ldots, y_{a_n}] \otimes \kk[R]$. Our choice of $a_1, \ldots, a_n$ implies that $\ZZ \ee_{a_1} \oplus \cdots \oplus \ZZ\ee_{a_n}$ maps injectively into $\ZZ^n$ and so $\kk[\NN^{Q_1} \oplus R] / I_\LL \cong \kk[y_{a_1}, \ldots, y_{a_n}] \otimes \kk[R]$. Therefore $[\VV(I_\LL) / \PGL(\alpha)] \cong [\AA^n \times (\kk^\times)^R / \PGL(\alpha)]$. 

We note that we may always fix the $(\kk^\times)^R$ component to 1. Now, the characters of the subgroup of $\PGL(\alpha)$ fixing the $(\kk^\times)^R$ component are given by $\Wt(Q) / R$. The map $\pic$ is surjective onto $G^\vee$ and its kernel is given by $R$, so  that $\Wt(Q) / R \cong G^\vee$. Hence the  aforementioned subgroup is naturally isomorphic to $G$. Consequently, we have stack isomorphisms \[ [\VV(I_\LL) / \PGL(\alpha)] \cong [\AA^n \times (\kk^\times)^R / \PGL(\alpha)] \cong [\AA^n \times \{1\} /G] \cong [\AA^n /G].\] This completes the proof.
\end{proof}

We recall the definitions of $G$-Hilb$(\AA^n)$ and Hilb$^G(\AA^n)$. Take $n\in \NN$ and $G$ as above. Following Reid \cite{GHilb}, define $G$-Hilb$(\AA^n)$ to be the fine moduli space of $G$-invariant subschemes of $\AA^n$ whose coordinate ring is isomorphic to $\kk[G]$ as a $\kk[G]$-module. Although the scheme $G$-Hilb($\AA^n$) is reducible in general, it has a distinguished irreducible component Hilb$^G(\AA^n)$ birational to $\AA^n/G$, see Ito-Nakumra \cite{ItoNakamura}. When $n \leq 3$ and $G \subset \SL(n,\kk)$, the scheme $G$-Hilb$(\AA^n)$ is smooth and isomorphic to Hilb$^G(\AA^n)$; furthermore the map $\tau: G\text{-Hilb}(\AA^n) \rightarrow \AA^n/G$ sending a subscheme to the orbit supporting it, is a crepant resolution of $\AA^n/G$, see \cite{ItoNakamura} and Nakamura \cite{Nakamura}.

Proposition \ref{clsdimm} allows us to recover the stack $[\AA^n / G]$ from the McKay quiver. Craw-Maclagan-Thomas \cite{CMT1} show that the distinguished component Hilb$^G(\AA^n)$ of $G$-Hilb$(\AA^n)$ can also be recovered from the labelled McKay quiver. Indeed, by Proposition 5.2 of \cite{CMT1}, Hilb$^G(\AA^n)$ is the subvariety of \[\MM_\vartheta(Q) = \big(\AA^{Q_1}\big) ^{ss} _\vartheta\, /\, \PGL(\alpha)\] cut out by the ideal \[I_Q:= \Big\langle y^{u_1} - y^{u_2} \, \Big|\, \div(u_1 - u_2) = 0, \,\inc(u_1 -u_2) = 0 \Big\rangle.\]   We proceed to relate our construction to that of Hilb$^G(\AA^n)$ by defining a GIT problem in which $[\AA^n /G]$ and Hilb$^G(\AA^n)$ are separated by a finite series of wall-crossings. We begin by carefully picking a basis $\BB$ of $R$.

Write $G^\vee$ as a direct sum of cyclic groups $\bigoplus_{j=1}^mH_j$ and take $\rho_j$ a generator of $H_j$. Define
\begin{equation*}
\overline{\BB}:= \Big\{-\ee_{\rho_j} -\ee_{\rho'\otimes \rho_j^{-1}} + \ee_{\rho'} \in \Wt(Q) \,|\, \forall\,1\leq j\leq m, \,\, \rho' \in G^\vee \setminus \{\rho_1, \ldots, \rho_m\}\Big\}.
\end{equation*}

\begin{lemma}
The set $\overline{\BB}$ generates the lattice $R \subset \Wt(Q)$.
\end{lemma}

\begin{proof}
For notational purposes, we use + for the binary operation on $G^\vee$ in this proof. First we show that
\begin{equation*}
\widetilde{\BB}:= \Big\{-\ee_{\rho_j} -\ee_{\rho' - \rho_j} + \ee_{\rho'} \in \Wt(Q) \,|\,\, \forall\,1\leq j\leq m,\,\,\, \rho' \in G^\vee\Big\}
\end{equation*}
generates $R$. Let $\rho = \sum_{j} \gamma_j \rho_j$ and without loss of generality assume $\gamma_j >0$. Since
\begin{equation*}
\sum_{1 \leq \kappa_j \leq \gamma_j} -\ee_{\rho_j} - \ee_{\rho' - (\kappa_j -1)\rho_j -\rho_j} +\ee_{\rho' - (\kappa_j -1)\rho_j} = -\gamma_j \ee_{\rho_j} - \ee_{(\rho' - \gamma_j \rho_j)} +\ee_{\rho'} 
\end{equation*}
we that deduce $(\sum_j -\gamma_j \ee_{\rho_j}) + \ee_{\rho'}$ is an element of $\NN \widetilde{\BB}$. Moreover, for $\rho' = \sum_{j} \gamma_j \rho_j$ and $\rho'' = \sum_{j} \gamma_j' \rho_j$, we have that $-\ee_{\rho'} - \ee_{\rho''} + \ee_{\rho' +\rho''}$ is equal to
\begin{equation*}
\Big(\Big(\sum_j \gamma_j \ee_{\rho_j}\Big) - \ee_{\rho'}\Big)  +\Big(\Big(\sum_j \gamma_j' \ee_{\rho_j}\Big) - \ee_{\rho''}\Big) + \Big(\Big(\sum_j -(\gamma_j+\gamma_j') \ee_{\rho_j}\Big) + \ee_{\rho' + \rho''}\Big)
\end{equation*}
showing that $-\ee_\rho' - \ee_{\rho''} + \ee_{\rho'+\rho''} \in \ZZ \widetilde{\BB}$. Therefore $\widetilde{\BB}$ generates $\ker(\div) =R$.

Take $-\ee_{\rho_j} - \ee_{\rho_j' - \rho_j} + \ee_{\rho_j'}$ for $1 \leq j, j' \leq m$ and $|\rho_j|$ to be the order of $\rho_j$. Then we have
\begin{equation*}
\sum_{0\leq \kappa \leq |\rho_j|- 2} -\ee_{\rho_j} - \ee_{\kappa \rho_j + \rho_{j'}} +\ee_{(\kappa+1)\rho_j + \rho_{j'}} = (1- |\rho_j|) \ee_{\rho_j} - \ee_{\rho_{j'}} + \ee_{\rho_{j'} - \rho_j}
\end{equation*}
which along with the fact that $|\rho_j|\, \ee_{\rho_j} \in \ZZ \overline{\BB}$ shows that $\overline{\BB}$ generates $R$.
\end{proof}

\begin{rmk}\label{P}
Note that for any $j$ as above, the $\ee_{\rho_j}$ coefficient of elements of $\NN \overline{\BB}$ is non-positive. This will prove crucial in the proof of the theorem below. 
\end{rmk}

Fix a basis $\BB \subset \overline{\BB}$ of $R$. We have the following diagram
\begin{equation}\label{com6}\begin{split}
\xymatrix@C=1.5cm@R=1cm{\NN^{Q_1} \oplus \NN \BB \ar[r]^-{\inc \oplus \iota} \ar[d]_{\div} & \Wt(Q) \\ \NN^n}
\end{split}\end{equation}
The semigroup homomorphism induces a $\Wt(Q)$-grading on $\kk[\NN^n \oplus \NN \BB]$ and hence an action of $\PGL(\alpha)$ on $\AA^{Q_1} \times \AA^\BB$. Define $I_{\LL, \BB}$ to be the $\Wt(Q)$-homogenous ideal of $\kk[\NN^n \oplus \NN \BB]$ given by
\begin{equation*}
I_{\LL, \BB} := \Big\langle y^{u_1}z^{v_1} - y^{u_2}z^{v_2} \, \Big| \, \div(u_1-u_2) = 0, \inc(u_1-u_2)+\iota(v_1-v_2)=0 \Big\rangle.
\end{equation*}
For $\theta \in \Wt(Q)$ we consider the stack quotient $[\VV(I_{\LL, \BB}) ^{ss} _\theta / \PGL(\alpha)]$.

\begin{thm}
There exists generic stability conditions $\chi_{\theta_1}, \chi_{\theta_2} \in \PGL(\alpha)^\vee$, such that \[[\AA^n/G] \cong [\VV(I_{\LL, \BB})^{ss}_{\theta_1} / \PGL(\alpha)] \,\text{ and }\, \textup{Hilb}^G(\AA^n) \cong [\VV(I_{\LL, \BB})^{ss}_{\theta_2} / \PGL(\alpha)].\]
\end{thm}

\begin{proof}
Because of Proposition \ref{clsdimm}, to establish the first isomorphism it suffices to find $\theta_1$ for which $\VV(I_{\LL , \BB})^{ss}_{\theta_1} = \VV(I_{\LL})_{\theta_1}^{ss}$. The cone $\QQ_{\geq 0} \BB \subset \Wt(Q)_\QQ$ is top dimensional, so we may pick a generic $\theta_1 \in \NN \BB$. After picking a higher multiple if necessary, we may assume that the $\chi_\theta$-unstable locus in $\AA^{Q_1} \times \AA^\BB$ is given by the vanishing locus of the ideal
\begin{equation*}
B_{\theta_1}:= \Big\langle\, y^uz^v \in \kk[\NN^{Q_1} \oplus \NN \BB] \, \Big| \, \inc(u) + \iota(v) = \theta_1 \,\Big\rangle
\end{equation*}
We claim that for any monomial $y^uz^v \in B_{\theta_1}$ there exists $u' \in \NN^{Q_1}$ such that $y^{u'}z^{\theta_1} - y^uz^v \in I_{\LL, \BB}$. Since $\theta_1 \in \NN \BB$, $\inc(u) = \theta_1 - \iota(v) \in R$. The commutative diagrams (\ref{com5}) imply $\div(u)$ is a torus-invariant section of the trivial line bundle. Take $u' \in \NN^{Q_1} $ to be a cycle in the quiver of sections corresponding to $\div (u)$ and note that $\inc(u')=0$. We then have that $\div(u-u') = 0$ and $\inc(u-u') + \iota(v-\theta_1) =0$, as claimed. From this it follows that any point for which $z^{\theta_1} =0$ is unstable. Now, $\theta_1$ is in the interior of $\QQ_{\geq 0}\BB$ therefore any point for which $z_\b = 0$ is unstable. That is  $\VV(I_{\LL, \BB})^{ss}_{\theta_1} := \VV(I_{\LL, \BB}) \setminus \VV(B_{\theta_1}) \subset \AA^{Q_1} \times (\kk^\times)^\BB$ and hence \[\VV(I_{\LL, \BB})^{ss}_{\theta_1} = \VV(I_\LL)^{ss}_{\theta_1}.\]

Now take a generic $\theta_2 \in \Wt(Q)$ in the interior top-dimensional cone $\Theta:= \QQ_{\geq 0}\{\ee_\rho \,|\, \rho \in G^\vee\setminus\{0\}\}$. Once again, taking a higher multiple if necessary we may assume that the $\chi_\theta$-unstable locus in $\AA^{Q_1} \times \AA^P$ is given by the vanishing locus of the ideal 
\begin{equation*}
B_{\theta_2}:= \Big\langle\, y^uz^v \in \kk[\NN^{Q_1} \oplus \NN P] \, | \, \inc(u) + \iota(v) = \theta_2 \,\Big\rangle.
\end{equation*}
If we set
\begin{equation*}
B_{\theta_2}':= \langle\, y^u \in \kk[\NN^{Q_1}] \, | \, \inc(u) = \theta_2\,\rangle
\end{equation*}
then the vanishing locus of $B_{\theta_2}'$ is equal to that of $B_\vartheta$, since $\theta_2$ and $\vartheta$ lie in the same chamber. Let $y^uz^v \in B_{\theta_2}$ and take $y^{u'}$ to be the unique monomial in $\kk[\NN^{Q_1}] / I_Q$ for which $\div(u) =\div(u')$ and $\inc(u') = \theta_2$. We next show that
\begin{equation}\label{iso}
\bigg(\frac{\kk[\NN^{Q_1} \oplus \NN B]}{I_{\LL, B}}\bigg)_{y^uz^v} \cong \bigg(\frac{\kk[\NN^{Q_1}]}{I_Q}\bigg)_{{y^{u'}}}.
\end{equation}
Remark \ref{P} gives that $v$ has a non-positive coefficient for each basis element $\ee_{\rho_j}$. Since $\inc(u) = \theta_2 - \iota(v)$ and $\theta_2$ is in the interior of $\Theta$, $\inc(u)$ has a strictly positive coefficients for each basis element $\ee_{\rho_j}$. Write $u = u_1 + \cdots u_m + u''$ for $u_j, u'' \in \NN^{Q_1}$ satisfying $\inc(u_j) = \ee_{\rho_j}$. We have that $y^{u} = y^{u_1}\cdots y^{u_m} y^{u''}$ and therefore in the localization above the monomials $y^{u_j}$ are invertible. Take an arbitrary element $\b:= -\ee_{\rho_j} - \ee_{\rho'} + \ee_{\rho_j \otimes \rho'}$ of $\BB$. Then there exists a path $p_j$ from $\rho'$ to $\rho' \otimes \rho_j$  with label $\div(u_j) \in \Hom(\rho', \rho'\otimes\rho_j)$. Let $u_j' \in \NN^{Q_1}$ be the element determined by $p_j$. We then have $\div(u_j - u_j')=0$ and $\inc(u_j - u_j') + \iota(\b) =0$ which implies that $z_\b y^{u_j}- y^{u_j'} \in I_{\LL,\BB}$. Since $y^{u_j}$ is invertible in the localization we may replace $z_\b y^{u_j}- y^{u_j'}$ by $z_\b- y^{u_j'-u_j}$, thereby eliminating $z_\b$ for every $\b\in \BB$. Next, consider the general generator $y^{u_1}z^{v_1} - y^{u_2}z^{v_2}$ of $I_{L,\BB}$, eliminate the monomials $z^{v_1}, z^{v_2}$ and multiply by the invertible elements $y^{u_j}$ to get a polynomial. Since every $z_\b$ is replaced by some $y^{u_j'-u_j}$ for which $\div(u_j - u_j')=0$ the resulting polynomial will be in $I_Q$. After noting that the construction above enables us to write $y^uz^v$ as a monomial in $\kk[\NN^{Q_1}]$ we have the isomorphisms (\ref{iso}).

The isomorphisms (\ref{iso}) allow us to conclude that  \[\VV(I_{\LL, \BB})^{ss}_{\theta_2} \cong \VV(I_Q)_{ss}^{\theta_2}.\] Consequently $[\VV(I_{\LL, \BB})_{\theta_2}^{ss} /  \PGL(\alpha)] \cong [\VV(I_Q)_{\theta_2}^{ss} / \PGL(\alpha)]$. Now  $[\VV(I_Q)_{\theta_2}^{ss} / \PGL(\alpha)]$ is a substack of the representable stack $\MM_{\theta_2}(Q)$ cut out by the homogenous ideal $I_Q$ and is therefore the variety $\VV(I_Q)_{\theta_2}^{ss} / \PGL(\alpha)$. Proposition 5.2 of \cite{CMT1} then gives Hilb$^G \cong [\VV(I_Q)_{\theta_2}^{ss} / \PGL(\alpha)]$, completing the proof. \end{proof}

\begin{coro}
For $n \leq 3$ and $G\subset \SL(n, \kk)$ finite abelian, there exists generic stability conditions $\chi_{\theta_1}, \chi_{\theta_2} \in \PGL(\alpha)^\vee$, such that \[ [\AA^n/G] \cong [\VV(I_{\LL, \BB})^{ss}_{\theta_1} / \PGL(\alpha)] \,\text{ and }\, G\textup{-Hilb}(\AA^n) \cong [\VV(I_{\LL, \BB})^{ss}_{\theta_2} / \PGL(\alpha)].\]
\end{coro}
\begin{proof}
For  n=2, by the work of Ito-Nakamura \cite{ItoNakamura} we have $G\textup{-Hilb}(\AA^n) \cong \textup{Hilb}^G(\AA^n)$. For n=3, Nakamura \cite{Nakamura} shows $G\textup{-Hilb}(\AA^n) \cong \textup{Hilb}^G(\AA^n)$ 
\end{proof}

\begin{rmk}
The careful choice of basis $\BB$ was made with $G\textup{-Hilb}(\AA^n)$ in mind. Moving between $[\AA^n/G]$ and another crepant resolution via wall-crossings may require a different choice of basis.
\end{rmk}

\addcontentsline {toc} {section} {Bibliography}
\bibliographystyle {plain}
\bibliography{Library}
\end{document}